\documentclass[11pt]{amsart}
\usepackage{amsmath,amssymb,amsfonts,epsf,epsfig}
\usepackage[labelfont=normal]{subcaption}
\usepackage{enumerate,amsthm,color,mathdots}
\usepackage{euscript,graphicx}
\usepackage{url, float, hyperref, color} 
\usepackage[usenames,dvipsnames,svgnames,table]{xcolor}
\usepackage{graphicx, mathtools}
\usepackage{tikz}

\newtheorem{theorem}{Theorem}[section]
\newtheorem{lemma}[theorem]{Lemma}

\newtheorem{corollary}[theorem]{Corollary}

\newtheorem{definition}[theorem]{Definition}

\newtheorem{remark}[theorem]{Remark}

\newcommand{\rev}[1] {\overline{#1}}

\newcommand{\Aut}{{\rm Aut}}

\newcommand{\GCD}{\mathrm{gcd}}

\newcommand{\Mon}{\mathrm{Mon}}

\newcommand{\ZZ}{\mathbb{Z}}

\newcommand{\F}{{\mathcal{F}}}
\newcommand{\M}{{\mathcal{M}}}
\newcommand{\N}{{\mathcal{N}}}
\newcommand{\W}{{\mathcal{W}}}
\newcommand{\Pe}{\hbox{\rm P}}

\newcommand{\Du}{\hbox{\rm D}}

\newcommand{\sk}{\hbox{\rm Sk}}

\title{Orbits of consistent  walks  in dart-transitive maps}
\author{Micael Toledo}
\author{Alejandra Ramos-Rivera}
\author{Primo\v{z} Poto\v{c}nik}
\author{Stephen E. Wilson}

\address{Primo\v{z} Poto\v{c}nik, Faculty of Mathematics and Physics, University of Ljubljana, Slovenia
\newline
\indent Also affiliated with: Institute of Mathematics, Physics and Mechanics, Ljubljana, Slovenia
}
\email{primoz.potocnik@fmf.uni-lj.si}

\address{Alejandra Ramos-Rivera, Institute of Mathematics, Physics and Mechanics, Ljubljana, Slovenia}
\email{alejandra.ramosrivera@fmf.uni-lj.si}

\address{Micael Toledo, Faculty of Mathematics and Physics, University of Ljubljana, Slovenia
\newline
\indent Also affiliated with: Institute of Mathematics, Physics and Mechanics, Ljubljana, Slovenia
}
\email{micaelalexitoledo@gmail.com}

\address{Stephen E. Wilson, Northern Arizona University, Flagstaff, Arizona, USA}
\email{stephen.wilson@nau.edu}

\thanks{The first named author was supported by ... The second- and third-named author acknowledges the support of the Slovenian Research Agency, programme number P1-0294 and project number J1-4351. The fourth named author wishes to thank for the hospitality and support the University of Ljubljana.}

\begin{document}
\maketitle


\begin{abstract}
In a simple graph, a {\em shunt} is a symmetry which sends an edge to an incident edge (without fixing their shared vertex).  The orbit of this edge under the shunt forms a {\em consistent cycle}.   The important theorem of Biggs and Conway says that in a dart-transitive graph of valence $q$, there are exactly $q-1$ orbits of consistent cycles.  These ideas  have become a useful tool in the area of graphs symmetries, and generalize easily to consistent walks in  graphs which are not simple.  These walks are not necessarily cycles, or even circuits.  
 This paper considers  these walks and their orbits in  the venue of dart-transitive maps and classifies them geometrically.
 \end{abstract}
 

\section{Introduction}
\label{sec:intro}


The starting point of this paper is the concept of a consistent cycle, which was first introduced by Biggs in \cite{B}.  In that paper,  Biggs states a theorem which he attributes to John Horton Conway \cite{conway}.  Its loose statement is that if a graph is simple, finite, symmetric and $q$-valent, then its symmetry group has exactly $q-1$ orbits of consistent cycles. Let us be more exact:

	Let $\Gamma$ be a finite simple graph. If $u,v$ are two adjacent vertices in $\Gamma$, then the pair $(u,v)$ is called a {\em dart} (or an {\em arc}) of $\Gamma$. 
 A {\em symmetry} or {\em automorphism} of $\Gamma$ is a permutation $\sigma$  of its vertices which preserves its edges. The set $\Aut(\Gamma)$  of all symmetries of $\Gamma$ is a group under composition.   It is important to realize that in the Biggs-Conway theorem, the word ‘cycle' refers to what we might call a {\em directed} cycle.  For clarity, we will use the word {\em ‘cyclet’} and define a cyclet to be a {\em sequence} $\alpha = [v_0, v_1, v_2, \dots, v_{r-1}]$, $r\ge 3$, of distinct vertices  such that for all $i\in \ZZ_r$, the vertex $v_i$ is adjacent to $v_{i+1}$ (with indices computed modulo $r$).  The cyclet is said to be {\em consistent} provided that for some symmetry $\sigma$ (called a {\em shunt} for $\alpha$), $v_i^\sigma = v_{i+1}$ for all $i\in\ZZ_r$. For a subgroup $G$ of $\Aut(\Gamma$) we say that $\alpha$ is $G$-consistent if it has a shunt in $G$.
The Biggs-Conway theorem then states the following:

\begin{theorem}[\cite{B,conway}] \label{th:BC1}
If $\Gamma$ is a finite simple $q$-valent graph and $G$ a subgroup of $\Aut(\Gamma)$ acting transitively on the darts  of $\Gamma$, then there are exactly $q-1$ orbits under $G$ of $G$-consistent cyclets.
\end{theorem}

While the theorem was largely unremarked at first, in the last two decades, interest in it  has increased, with new proofs being found, new variations introduced, and new uses constructed:   see \cite{BMP, GKM, K, KKJ, KM, MPW1, MPW2} for examples.
\smallskip

Note the Biggs-Conway theorem determines the number of orbits of consistent cycles but it gives no additional information about them (for example, what are the lengths of the consistent cycles or how can one find them in the graph).
Our intention here is to address this problem in a specific venue of (skeletons) of symmetric maps. In particular,
we will describe the orbits of consistent cycles in any dart-transitive map in terms of the $j$-hole-walks and $j$-Petrie-walks in the map.






Observe first that the skeleton of a map need not be a simple graph. Hence, if we are to investigate
consistent cycles in maps, we need to allow a more general notion of graphs, such as {\em multigraphs} (allowing multiple edges) or even {\em pseudographs} (allowing multiple edges and loops). We must recognize that in these venues,  the orbit of a shunt need not be a cycle, but might be something more general.

To cope with this, define a {\em dart} or {\em arc} to be one of the two ways to indicate a direction along an edge from one end to the other.  A dart has an initial vertex and a terminal vertex (and if an edge is a loop, those are the same vertex).  A pair of darts $(d_1, d_2)$ is {\em sequential} provided that they are not reverses of each other and the terminal vertex of $d_1$ is the initial of $d_2$.   We can now re-define the word ‘cyclet’:  a {\em cyclet} is a sequence $\alpha = [d_1, d_2, \dots, d_r]$ of darts such that for all $i\in\ZZ_r$, $(d_i, d_{i+1})$ is sequential. The cyclet is said to be {\em consistent} provided that for some symmetry $\sigma$ (called a {\em shunt} for $\alpha$), $d_i^\sigma = d_{i+1}$ for all $i\in\ZZ_r$.

With this definition of a cyclet, the Biggs-Conway theorem easily generalizes to pseudographs (the proof of the following generalization, which follows closely the proof of the Biggs-Conway theorem as presented in \cite{MPW2},  can be found in the note \cite{MPW3}):
\begin{theorem}[\cite{MPW3}] \label{th:BC2}
{\em If $\Gamma$ is a finite pseudograph, $G\le \Aut(\Gamma)$ acts transitively on the darts  of $\Gamma$ and $q$ is the common valence of vertices in $\Gamma$, then there are exactly $q-1$ orbits under $G$ of $G$-consistent cyclets.}
\end{theorem}
With this generalization in mind we will give more information on how (orbits of) consistent cycles emerge in (the skeletons of) dart-transitive maps.



The paper is organized as follows: in Section~\ref{sec:maps}, we give a combinatorial definition of a map in terms of its flags and their connections.  We also provide a more geometric/topological definition of a map as an embedding of a pseudograph in a surface.  We discuss  symmetries (or automorphisms) of maps in both presentations.  In Section~\ref{sec:flags}, we define notions related to those of walks, circuits and cycles in graph theory.  We introduce Coxeter's ideas of $j$-th order ``holes" and ``Petrie paths''.  In Section~\ref{sec:DTMaps}, we talk about types of symmetry in a map and describe the five kinds of dart-transitive maps.  In  Section~\ref{sec:cw}, we intoduce the idea of {\em consistency} in closed walks and classify all that happen in maps.  We offer Theorem \ref{the:const}   which proves that every consistent walk is a $j$-hole or a $j$-Petrie path for some $j$.  With that established, we prove the first of our two important theorems, Theorem \ref{the:orbits}, which determines the number and kind of each orbit of consistent walks. In Section~\ref{sec:cbb}, we describe four shapes of  the edge-sets of consistent walks.  We prove our second important theorem, Theorem \ref{th:classholes}, proving that every consistent walk belongs to one of those four types.  In Section \ref{sec:onevertex},  we prove the corresponding results in  maps having exactly one vertex.
\color{black}

Let us end this introduction with a few items of notation.
We will use exponential notation for permutations: if $\sigma$ is an element of a permutation group $G$ acting on a set $X$ and $x\in X$, we write $x^\sigma$ for the image of $x$ under $\sigma$, and $x^G$ for the orbit of $x$ under $G$. Moreover, if $S = (x_1,x_2,\ldots,x_n)$ is a sequence of elements of $X$, we let $S^g$ be the sequence $(x_1^g,x_2^g,\ldots,x_n^g)$. In many situations it will be convenient to adopt the notation $|a|_n$, read, ``the order of $a$ mod $n$'', to mean the additive order of the element $a$ in $\ZZ_n$.  This is equal to $\frac{n}{\gcd(n,a)}$.

\section{Maps and their symmetries}
\label{sec:maps}

Often a map is defined as an embedding of a graph (not necessarily simple) into a connected compact $2$-manifold without boundary, sometimes satisfying some additional requirements. The barycentric subdivison of such an embedding decomposes the surface into triangles, called flags. The connections between flags motivates the combinatorial  definition of a map in the next paragraphs. The reconstruction of the embedding from the corresponding combinatorial definition of a map is briefly explained in Subsection~\ref{sec:model}.

A {\em map}   is a 4-tuple $\M := (\F, r_0, r_1, r_2)$ where $\F$ is a non-empty set of {\em flags} and $(r_0,r_1,r_2)$ is a triple of fixed-point free permutations of the elements of $\F$ such that $r_0^2=r_1^2=r_2^2=(r_0r_2)^2=1$. We will additionally require that the action of the group $\Mon(\M) := \langle r_0,r_1,r_2\rangle$ is transitive and that the images of any given flag $\Phi \in \F$ under $r_0$, $r_1$ and $r_2$ are all distinct. For $i \in \{0,1,2\}$ we let $\Phi^i$ denote the image of a flag $\Phi$ under $r_i$, and we say that the flags $\Phi$ and $\Phi^i$ are {\em $i$-adjacent} or are {\em $i$-neighbours}. In the same spirit, we let $\Phi^{i_1i_2\ldots i_k}$ denote the image of $\Phi$ under the permutation $r_{i_1}r_{i_2}\ldots r_{i_k}$. For brevity, we may, for instance, write $\Phi^{(01)^3}$ instead of $\Phi^{010101}$. 

Observe that the group $\Mon^+(\M):=\langle r_0r_1,r_1r_2\rangle$ has index at most $2$ in $\Mon(\M)$; if it {\em is} $2$, then $\M$ is {\em orientable} 
and if $\Mon^+(\M)=\Mon(\M)$, then $\M$ is {\em non-orientable}. 

For $i \in \{0,1,2\}$, an $i$-face of $\M$ is an orbit of flags under the group $\langle r_j \mid j \in \{0,1,2\}\setminus \{i\} \rangle$. The $0$-, $1$- and $2$-faces of $\M$ are called vertices, edges and faces, respectively. For $i,j \in \{0,1,2\}$, $i\neq j$, an $i$-face is {\em incident} to a $j$-face if their intersection is non-empty. A {\em dart} is an orbit of flags under $\langle r_2 \rangle$. 

Let $x =\{\Phi,\Phi^2\}$ be a dart. We define the {\em initial vertex} of $x$ as the unique vertex $u$ containing both $\Phi$ and $\Phi^2$. The {\em reverse} of $x$ is the dart $\rev{x}=x^{r_0}=\{\Phi^0,\Phi^{20}\}=\{\Phi^0,\Phi^{02}\}$. Note that the ``reverse'' relation is symmetric: the reverse of the reverse of $x$ is $x$ itself. The {\em terminal vertex}  of $x$ is the initial vertex of $\overline{x}$. Observe that an edge consists precisely of the four flags contained in a dart $x$ and its reverse $\overline{x}$.

The {\em skeleton} of $\M$ is the graph $\textrm{Sk}(\M)$ whose vertices and darts are the vertices and darts of $\M$ and incidence is given by containment. The skeleton of a map is a connected graph that may admit parallel edges and loops.   Such an object is sometimes referred to as a {\em multigraph} or  a {\em pseudograph}. If $\sk(\M)$ has no parallel edges and no loops, we say that $\M$ is {\em simple}. 

The {\em valence} of a vertex is defined as the number of darts that it contains. This number equals half the number of flags that it contains.

\begin{definition}For any flag $\Phi$, let $e$ and $v$ be the edge and vertex to which it belongs, respectively.  Then for any $d$ less than the valence of $v$, the flags $\Phi^{(r_1r_2)^d}$ and $\Phi^{(r_1r_2)^{d-1}r_1}$ belong to one edge, call it $e'$, which is also incident with $v$.  In this case we say that $e$ and $e'$ {\em subtend} $d$ faces at $v$.    
\end{definition}

Of course, if $e$ and $e'$ subtend $d$ faces at $v$ and $v$ has valence $q$, they also subtend $q-d$ faces there, where $q$ is the valence of $v$.  We use this language even when there might be repetitions among the faces enclosed by $e$ and $e'$.

\subsection{Realizations}
\label{sec:model}
To provide some intuitive insight,
 we construct here a topological object which embodies and motivates the definition of an abstract map.

\begin{figure}[hhh]
\begin{center}
\epsfig{file=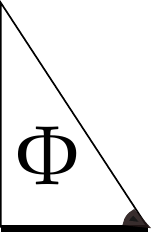,height=15mm}
\caption{The right triangle corresponding to the flag $\Phi$.}
\label{fig:OneFlag}
\end{center}
\end{figure}

Let $\M = (\F, r_0, r_1, r_2)$ be a map.   For each $\Phi\in\F$, create a right triangle (with interior) labelled $\Phi$, with one leg and the incident acute angle darkened, as in Figure~\ref{fig:OneFlag}.

\begin{figure}[hhh]
\begin{center}
\epsfig{file=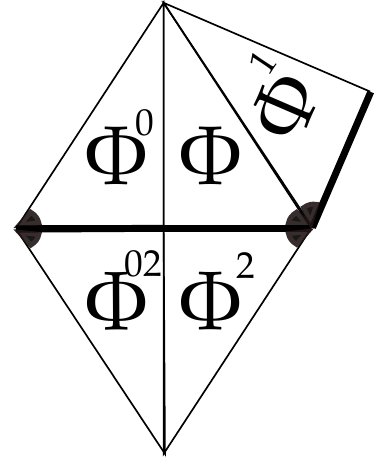,height=32mm}
\caption{The neighbours of the flag $\Phi$.}
\label{fig:AdjFlag}
\end{center}
\end{figure}

Then for each flag $\Phi$, attach  its triangle to those of $\Phi^0, \Phi^1, \Phi^2$ by identifying corresponding edges as shown in Figure~\ref{fig:AdjFlag}.

The resulting topological space is a 2-manifold, and the darkened parts are a topological form of a graph (possibly with loops and parallel edges).  This is a {\em realization} of $\M$, also called a {\em topological form} of $\M$. Note that $\M$ is orientable if and only if its
topological form is on an orientable surface.

\subsection{Automorphisms}

An {\em automorphism} or {\em symmetry}  of a map $\M$ is a permutation of the flags that commutes with the three involutions $r_0$, $r_1$ and $r_2$. 
Note that an automorphism of $\M$ has an induced action on the set of its darts that preserves the ``initial vertex'' and ``reverse'' operators. In this sense, every automorphism of $\M$ is also an automorphism of its skeleton. 
The automorphisms of $\M$ form a group $\Aut(\M)$ under (right) composition. 

An automorphism $g$ is a {\em reflexion} provided it maps some flag $\Phi$ to its $i$-neighbour for some $i \in \{0,1,2\}$, 
and it is a {\em rotation} if it maps some flag $\Phi$ to $\Phi^r$ where $r$ is an element of $\Mon^+(\M)$. The rotations of $\M$ form a subgroup $\Aut^+(\M)$ of $\Aut(\M)$ of index at most $2$, called the {\em rotation group} of $\M$. 

We say that a subgroup $G \leq \Aut(\M)$ is {\em reflexible} if for every flag $\Phi$ and every $i\in\{0,1,2\}$ there exists a reflection in $G$ mapping
$\Phi$ to $\Phi^i$.
Similarly, $G$ is {\em rotary} if for every flag $\Phi$ and all $i,j\in\{0,1,2\}$ there exists a symmetry in $G$ mapping $\Phi$ to $\Phi^{ij}$.
If $G$ is rotary but not reflexible, then it is {\em chiral}.
Further, we say $G$ is {\em face-reflexible} provided that every flag  can be mapped to its $0$- and $1$-neighbour by an element in $G$,
and is {\em half-reflexible} if it is face-reflexible but not reflexible.
In the case when $G=\Aut(\M)$, we say that $\M$ has the corresponding one of the above properties.



\section{Sequences of flags: holes and Petrie paths}
\label{sec:flags}
For a positive integer $n$, a {\em flag-walk} $W$ of length $n$ is a sequence of flags $(\Phi_0, \Phi_1, ,\ldots,\Phi_{n-1})$ such that for all $i \in \{0,\ldots,n-2\}$, the flags
 $(\Phi_i)^0$ and $\Phi_{i+1}$ belong to the same vertex. The flag-walk $W$ is {\em closed} if $(\Phi_{n-1})^0$ and $\Phi_0$
 are in the same vertex. 

For a flag-walk $W$, define its {\em partner} $W^2$ to be the sequence $(\Phi_0^2, \Phi_1^2, ,\ldots,\Phi_{n-1}^2)$.  Then  $W^2$ is a flag-walk if $W$ is, and $W^2$ is closed if $W$ is. 
Further, $W$ and $W^2$ induce the same sequence of darts.

\begin{figure}[h!]
\begin{center}
\begin{tabular}{ccc}
\includegraphics[width=0.3\textwidth]{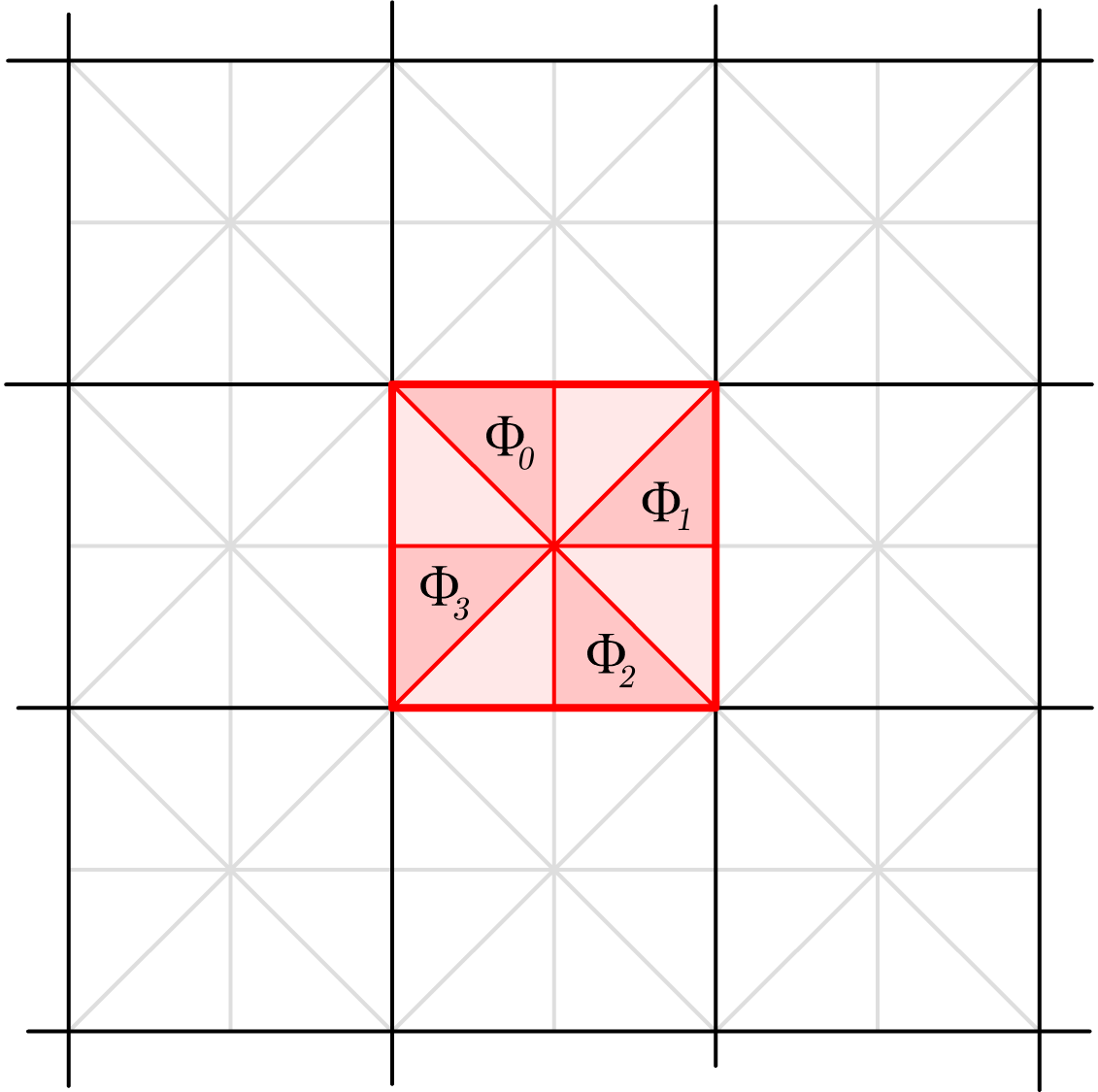} & 
\includegraphics[width=0.3\textwidth]{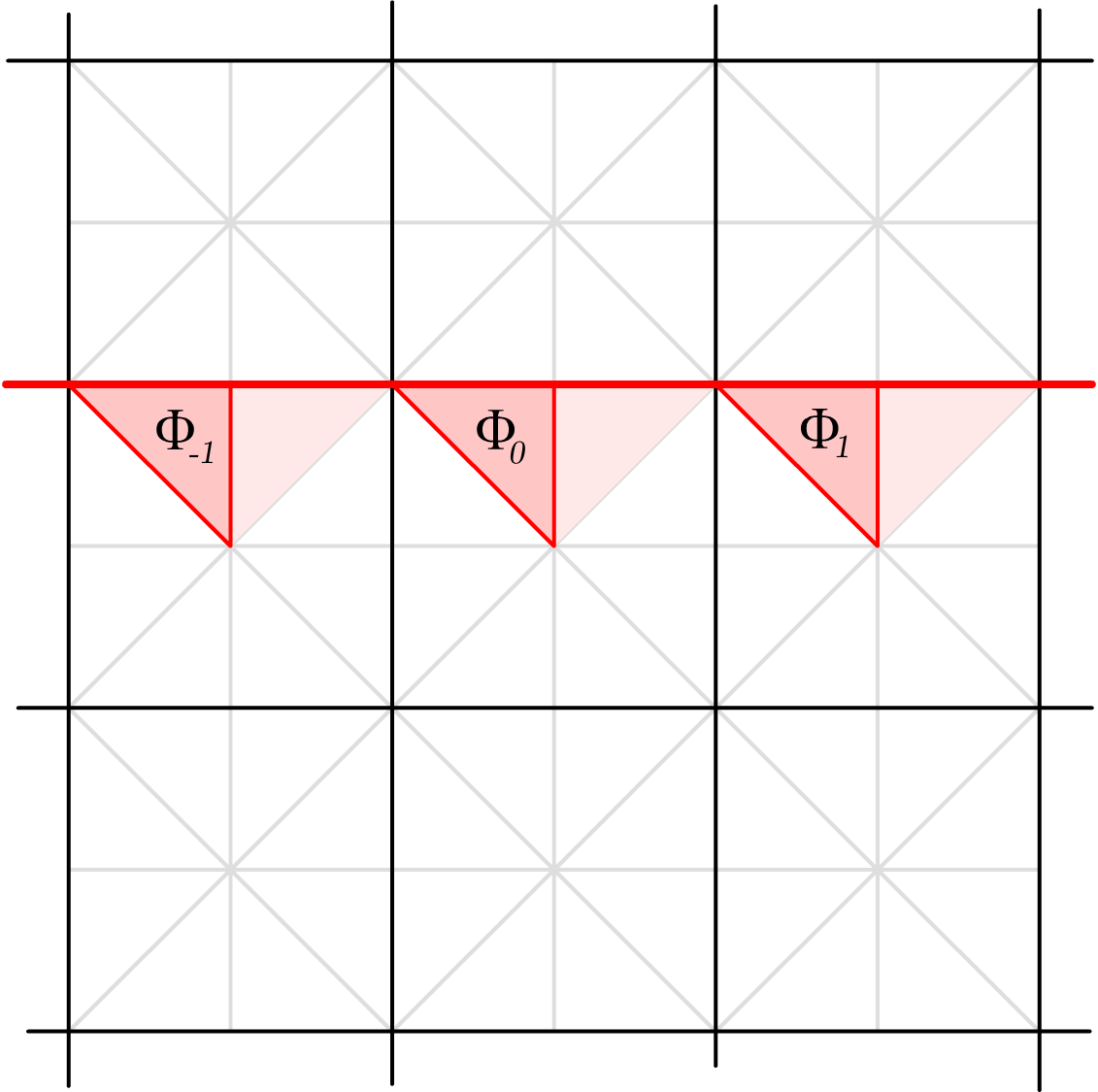} &  
\includegraphics[width=0.3\textwidth]{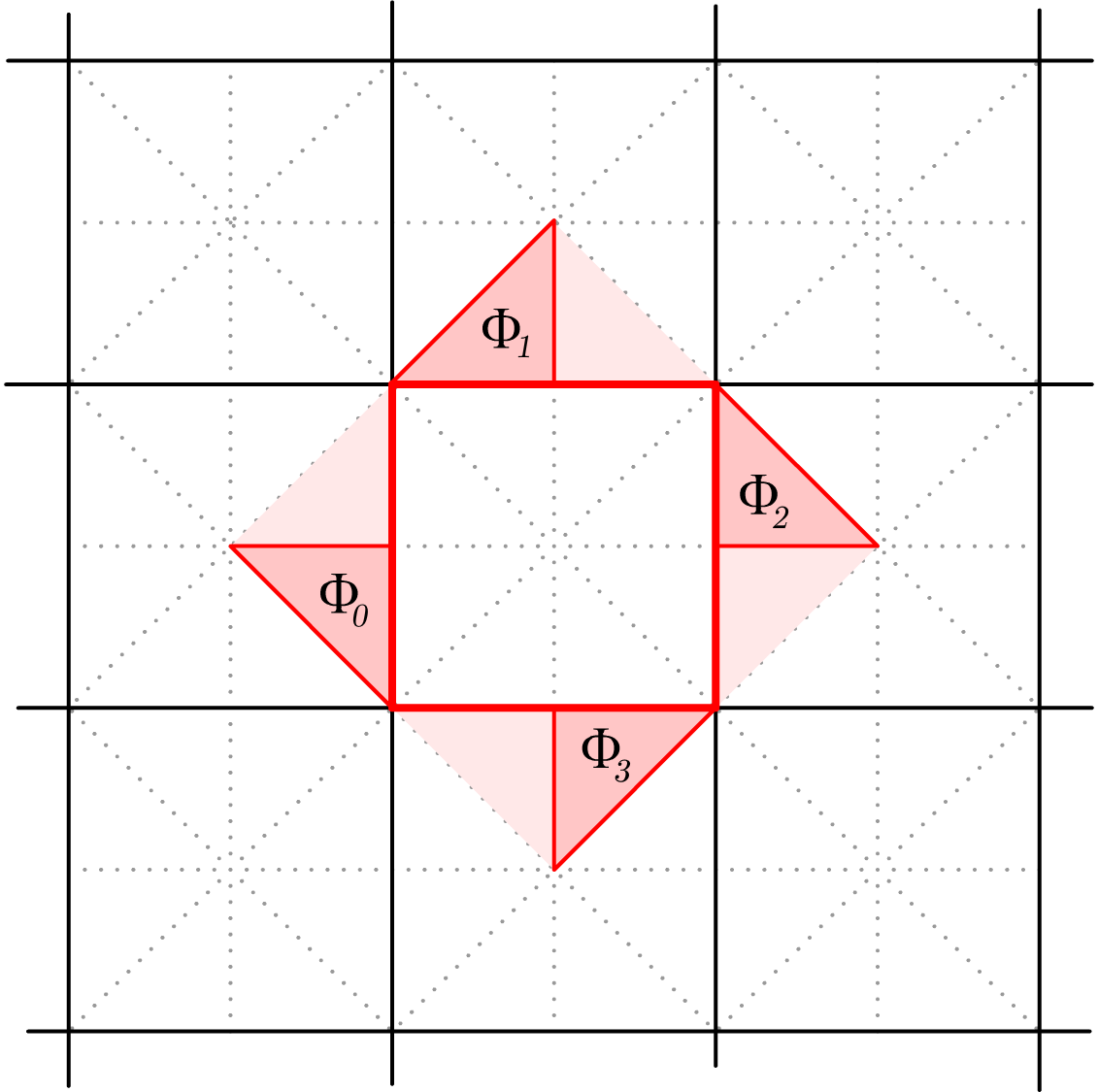}  \\
$1$-hole; $r=r_0r_1$ & $2$-hole; $r=r_0r_1r_2r_1$ & $3$-hole; $r=r_0r_1(r_2r_1)^2$ \\ 
\includegraphics[width=0.3\textwidth]{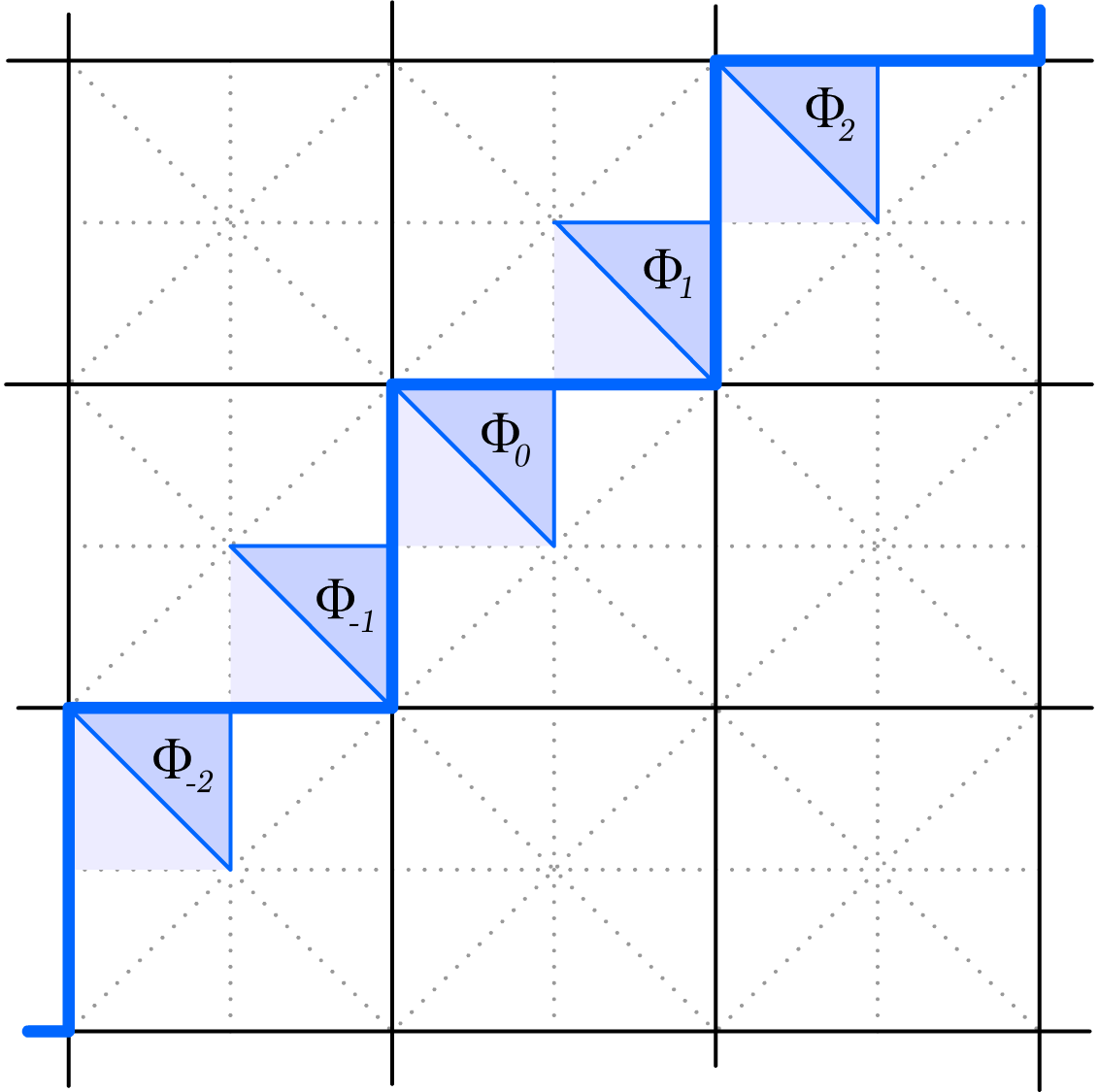} & 
\includegraphics[width=0.3\textwidth]{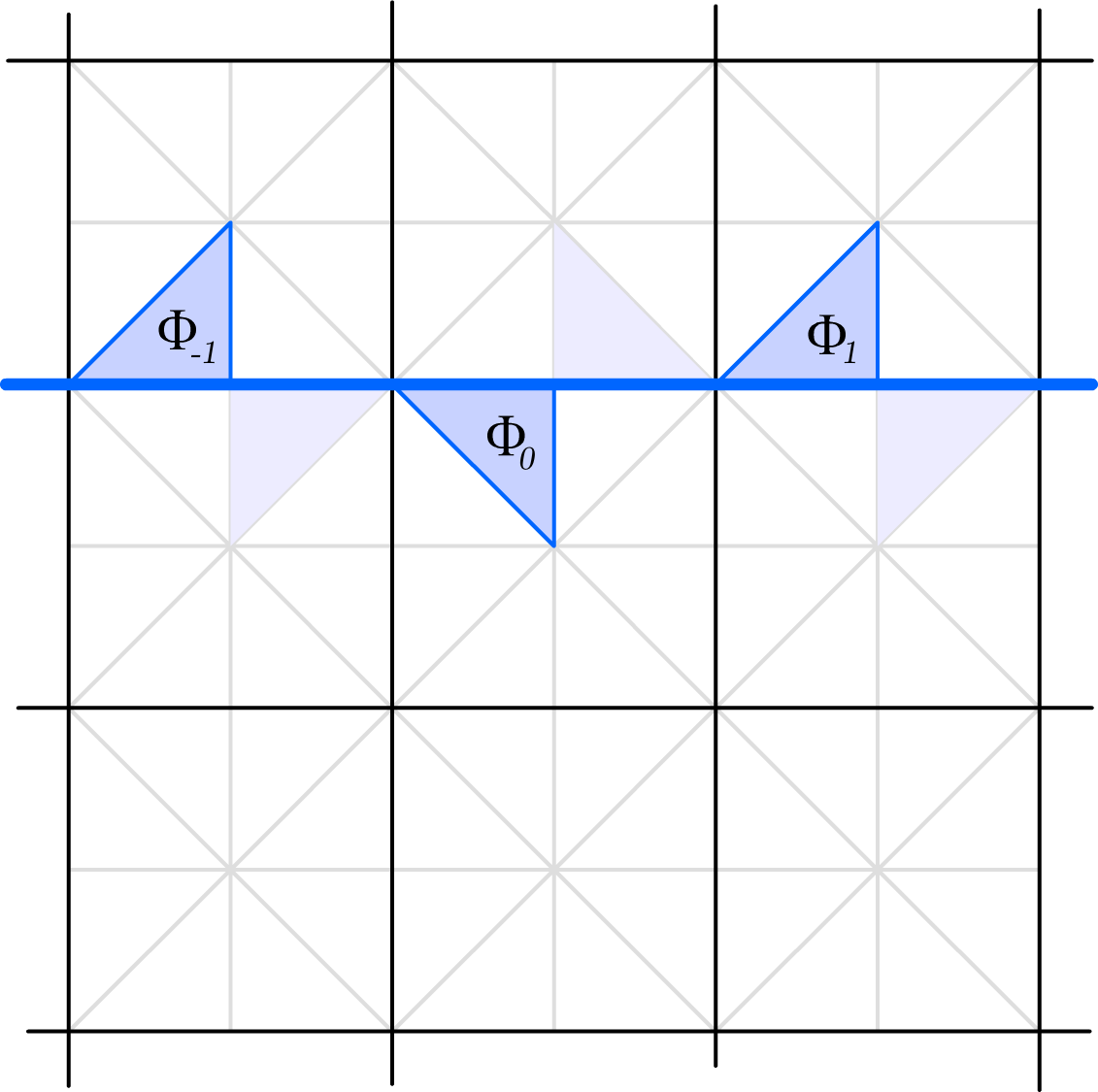} &  
\includegraphics[width=0.3\textwidth]{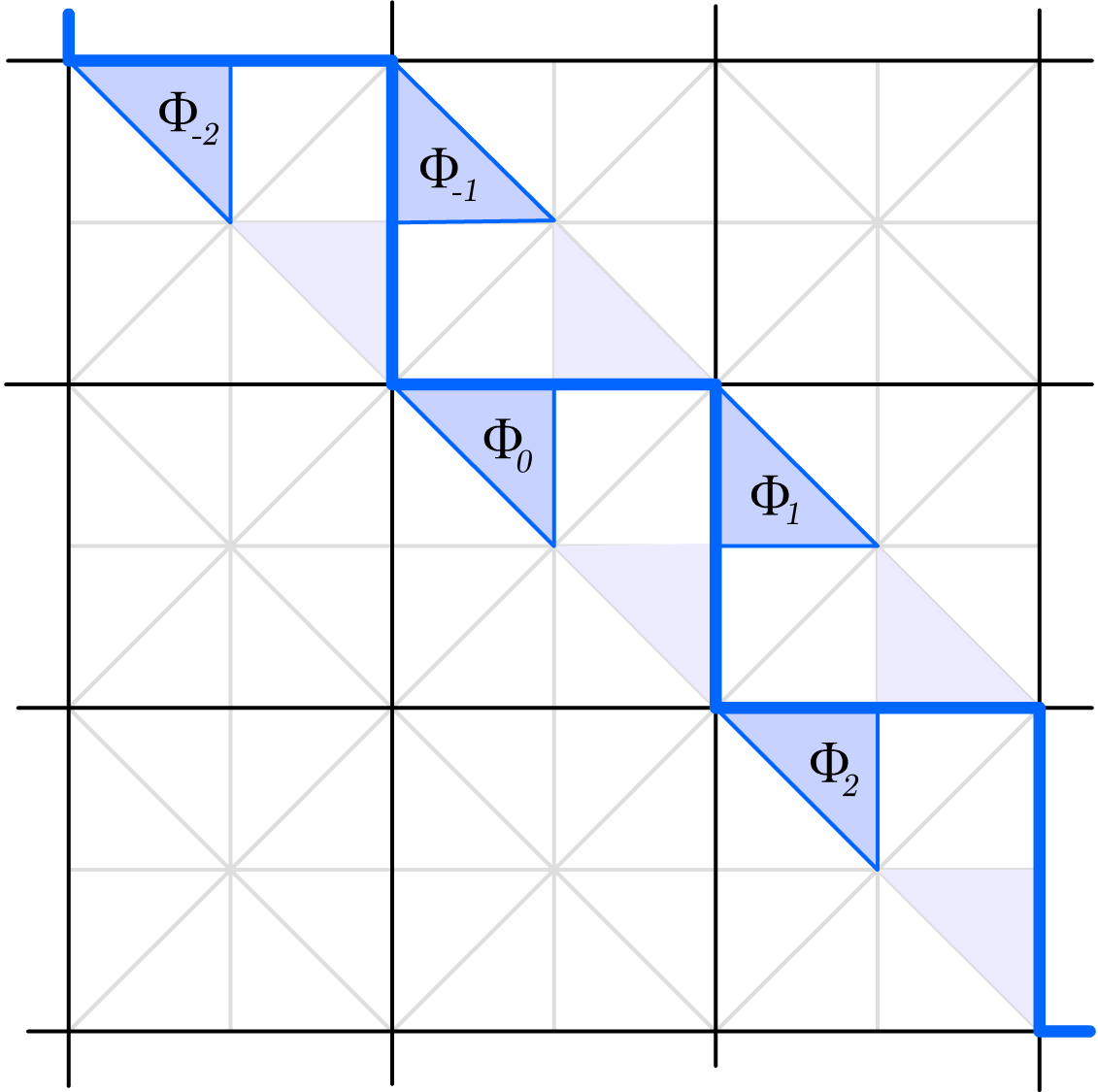}  \\
$1$-Petrie path; $s=r_0r_2r_1$ & $2$-Petrie path; $s=r_0(r_2r_1)^2$ & $3$-Petrie path; $s=r_0(r_2r_1)^3$
\end{tabular}
\caption{Different holes and Petrie paths on a $4$-valent map.}
\label{fig:holesandpetries}
\end{center}
\end{figure}

Assume that all vertices in a map $\M:=(\F,r_0,r_1,r_2)$ have valence $q$. 
For a positive integer $0<j < q$, let $\alpha_j$ be the permutation 
\begin{align}
\alpha_j = r_0r_1(r_2r_1)^{j-1}. 
\end{align}
In situations where $j$ is specified or understood, we may simply write $\alpha$.   Suppose that $\alpha$ has order $k$. For any flag $\Phi$, the sequence  $H = (\Phi_0, \Phi_1, \Phi_2, \dots, \Phi_{k-1})$, where $\Phi_0 = \Phi$ and $\Phi_{i+1} = \Phi^\alpha_i$, is  a {\em flag $j$-hole} (see top of Figure~\ref{fig:holesandpetries}).
 Compare this notation to that in Coxeter and Moser \cite{CM}.    A $j$-hole is always a closed flag-walk, and its partner $H^2$, then,  is a flag $(q-j)$-hole (compare the top left to the top right of Figure~\ref{fig:holesandpetries}). 

Further, the sequence $\rev{H} = (\Phi_{k-1}^0, \dots, \Phi_2^0, \Phi_1^0, \Phi_0^0)$ is called the {\em reverse} of $H$, and is a $j$-hole in its own right. The flags  $ \Phi_0, \Phi_1, \Phi_2, \dots, \Phi_{k-1}$,  together with the flags  $\Phi_0^0, \Phi_1^0, \Phi_2^0, \dots, \Phi_{k-1}^0$   are said to be the {\em interior} of the $j$-hole (see Figure~\ref{fig:holesandpetries}). That is, the interior of a $j$-hole $H$ consists of all the flags that are contained in either $H$ or its reverse $\rev{H}$. Note that the interior of a flag $1$-hole is a face of $\M$.

Now consider the permutation  $\beta_j$ (or just $\beta$) defined by 
\begin{align}
\beta_j= r_0(r_2r_1)^j
\end{align}
and suppose that it has order $\ell$.
For any flag $\Phi$,  let $P$ be the sequence $ (\Phi_0, \Phi_1, \Phi_2, \dots, \Phi_{\ell-1})$, where $\Phi_0 = \Phi$ and $\Phi_{i+1} = \Phi^\beta_i$. We will call it a {\em flag $j$-Petrie path}.  Note that the partner sequence $P^2 =  (\Phi_0^2, \Phi_1^2, \Phi_2^2, \dots, \Phi_{\ell-1}^2)$, is a flag $(q-j)$-Petrie path. Finally, the sequence $\rev{P}=(\Phi^{0 2}_{\ell-1}, \dots,  \Phi^{0 2}_2,  \Phi^{0 2}_1,\Phi^{0 2}_0) $, the {\em reverse} of $P$, is itself a $j$-Petrie path. The {\em interior} of $P$ is defined to be the set of all flags that belong to $P$ or $\rev{P}$.

In the special case where $\M$ is a map of even valence $q$, for any flag $\Phi$ the set of edges visited by a $\frac{q}{2}$-hole containing $\Phi$ is the same as that of the $\frac{q}{2}$-Petrie path containing $\Phi$. In this case, we call that set of a edges a {\em line} of $\M$.   

To get a more intuitive grasp of what $j$-holes and $j$-Petrie paths may look like, consider the following informal scenario. If we were to walk along the edges of a $j$-hole $H$, then all the flags belonging to $H$ would always appear on the same side (so either left or right) of the edges visited by $H$. Meanwhile, if we were to do the same on a Petrie path $P$, then flags would appear alternately on the left and  right side of the edges visited. In each case, consecutive edges of the $j$-hole or $j$-Petrie path subtend $j$ faces. This can be seen in Figure~\ref{fig:holesandpetries}, where the darkly coloured flags belong to a flag $j$-hole or flag $j$-Petrie path and the lightly coloured flags belong to the reverse of the corresponding hole or Petrie path. The dark coloured edges are the edges visited by the holes and Petrie paths. Observe that every flag $j$-hole and $j$-Petrie path visits the same set of edges as its reverse.

Let us now present two standard operations on maps, the {\em dual} and the {\em Petrie} (see, for example \cite{Wop}).
Both operations are involutory and produce a map from an existing map on the same set of flags. 

Let $\M = (\F,r_0,r_1,r_2)$ be a map. 
The {\em Petrie} of $\M$, denoted $\Pe(\M)$, is then the map $(\F,r_0r_2,r_1,r_2)$. 
Observe that each face of $\Pe(\M)$ is the interior of a flag $1$-Petrie path of $\M$. 
 
The dual of $\M$, denoted $\Du(\M)$, is the map $(\F,r_2,r_1,r_0)$. 
The faces of $\Du(\M)$ are the vertices of $\M$, and the vertices of $\Du(\M)$ are the faces of $\M$. 
  
Observe that $\Aut(\Du(\M))=\Aut(\M) = \Aut(\Pe(\M))$ and if one of $\M$, $\Pe(\M)$ or $\Du(\M)$ is reflexible, then they all are.

\begin{remark}
\label{rem:petrieholes}
The Petrie operator sends $j$-holes to $j$-Petrie paths and vice versa: the $j$-holes of $\M$ are the $j$-Petrie paths of $\Pe(\M)$ and the $j$-Petrie paths of $\M$ are the $j$-holes of $\Pe(\M)$. Similarly, the Petrie operator maps the interior of a $j$-hole to the interior of a Petrie path and vice versa.
\end{remark}

\color{Black}

\section{Dart-transitive maps}
\label{sec:DTMaps}

For a subgroup $G$ of $\Aut(\M)$, we say that a map $\M$ is {\em $G$-dart-transitive} provided that  $G $ acts transitively on the set of darts of $\M$. If $G = \Aut(\M)$ we may omit the prefix and simply say $\M$ is dart-transitive.



A map $\M$ is said to be a $k$-orbit map provided that the action of $\Aut(\M)$ on $\F(\M)$ has $k$ distinct orbits on flags.
In particular, reflexible maps are $1$-orbit maps and chiral maps are one kind of $2$-orbit maps.
For a subset $I \subsetneq \{0,1,2\}$, we say that a $2$-orbit map is in class $2_I$ given that two adjacent flags are in the same orbit if and only if they are $i$-neighbours for some $i \in I$. For convenience, instead of writing the set $I$ as a subindex, we will simply list its elements in increasing order. For instance, we will write $2$ instead of $2_{\emptyset}$, and $2_{01}$ instead of $2_{\{0,1\}}$.

It is known (see \cite{OPW}, for instance) that infinitely many maps exist in each of the seven possible classes of two-orbit maps. For $j \in \{0,1,2\}$, the two-orbit maps that are not $j$-face transitive are precisely those in class $2_I$ with $I = \{0,1,2\} \setminus \{j\}$; all others are transitive on the $j$-faces (recall that the $0$-, $1$- and $2$-faces of $\M$ are the vertices, edges and faces of $\M$ respectively).

\medskip

Let $\M$ be a $G$-dart-transitive map and let $x = \{\Phi, \Phi^2\}$ be a dart of $\M$. Let $\Psi$ be any flag of $\M$. By dart-transitivity, there exist a $g \in \Aut(\M)$ mapping $x$ to the dart $\{\Psi,\Psi^2\}$. In particular $\Phi$ must be in the same orbit as either $\Psi$ or $\Psi^2$. It follows that the action of $\Aut(\M)$ on $\F(\M)$ has at most two orbits. Observe that if $\Phi$ and $\Phi^2$ are in the same $\Aut(\M)$-orbit, then $\M$ is reflexible. Therefore, if $\M$ is not reflexible then no flag is in the same orbit as its $2$-neighbour, and thus $\M$ must be in class $2_I$ with $I \subseteq \{0,1\}$. That is, $M$ is either reflexible or in class $2$, $2_0$, $2_1$ or $2_{01}$.

The Petrie operator sends maps in class $2$ to maps in class $2_0$, and vice versa. Similarly, it sends maps in class $2_{1}$ to maps in class $2_{01}$ and vice versa.

Note that maps in class $2_0$ and $2_{1}$ are face-transitive but not reflexible, and the maps in class $2_{01}$ are half-reflexible. 

We summarize this discussion in Lemma~\ref{lem:darttrans}:


%

\begin{lemma}
\label{lem:darttrans}
A map $\M$ is dart-transitive if and only if one of the following holds:

\begin{enumerate}
\item $\M$ is reflexible; 
\item $\M$ is chiral (that is, $\M$ is in class $2$);
\item $\Pe(\M)$ is chiral (that is, $\M$ is in class $2_0$);
\item $\M$ is half-reflexible (that is, $\M$ is in class $2_{01}$)
\item $\Pe(\M)$ is half-reflexible (that is, $\M$ is in class $2_{1}$).
\end{enumerate}
\end{lemma}


We finish this section with the following remark, which follows easily from the definitions of different symmetry types of maps.

\begin{remark}\label{cor:refsub}
Suppose that $\M$ is a reflexible map.  Then:
\begin{itemize}
\item $\Aut(\M)$ has a chiral subgroup if and only if $\M$ is orientable.  That chiral subgroup is $\Aut^+(\M)$; it is the set of all orientation-preserving symmetries of $\M$.
\item $\Aut(\M)$ has a half-reflexible subgroup if and only if $\M$ is face-bipartite, meaning that its faces can be bi-coloured in such a way that adjacent faces have different colours. In this case, we call the subgroup $\Aut^\#(\M)$; this subgroup is the stabilizer in $\Aut(\M)$ of each chromatic class.
\end{itemize}
\end{remark}

\section{Consistent Walks}
\label{sec:cw}

Let $\M$ be a map, let $G \leq \Aut(\M)$ and let $W = (\Phi_0,\Phi_1,\ldots,\Phi_{n-1})$ be a closed flag-walk of length $n$. We say $W$ is a {\em $G$-consistent flag-walk} if there exists a $g \in G$ of order $n$ such that $\Phi_i^g = \Phi_{i+1}$ for all $i = 0, \ldots, n-2$. The automorphism $g$ is called a {\em shunt} of $W$. Note that since the action of $\Aut(\M)$ on the flags is semiregular and $g$ has order $n$, we have $\Phi_{n-1}^g = \Phi_0$.   If $W$ is a consistent flag-walk, the edge-set $\{e_0, e_1, \dots, e_{k-1}\}$ is simply a {\em consistent walk}, where each $\Phi_i$ is in edge $e_i$.

Our first task is to state and prove the following important theorem.

\begin{theorem}
\label{the:const}
Each consistent flag-walk in a map is, for some $j$, a flag $j$-hole or a flag $j$-Petrie path. 
\end{theorem}
\begin{proof} 

Let $W = (\Phi_0, \Phi_1, \Phi_2, \dots, \Phi_{k-1})$ be a consistent flag-walk in a map $\M$ and let $\sigma$ be a shunt for $W$. 
Then $\Phi_0^0$ and $\Phi_1$ are in the same vertex. 
Thus, for some $j$, $\Phi_1$ is $(\Phi_0^0)^{(2 1)^j}$ or $(\Phi_0^0)^{1(2 1)^{j-1}}$, and so $\Phi_1 = \Phi_0^{0(1 2)^j}$ or $\Phi_0^{01(2 1)^{j-1}}$.   
Thus $\Phi_1 = \Phi_0^\gamma$, where $\gamma$ stands for $\alpha_j=r_0r_1(r_2r_1)^{j-1}$ or $\beta_j=r_0(r_2r_1)^j$, as appropriate. 
But then $\Phi_2 = \Phi_1^\sigma = (\Phi_0^\gamma)^\sigma =(\Phi_0^\sigma )^\gamma = \Phi_1^\gamma$.  
Proceeding inductively, we see that $\W=(\Phi_0, \Phi_1, \dots, \Phi_{k-1})$ is a flag $j$-hole or a flag $j$-Petrie path. 
\end{proof}

There are two direct consequences of Theorem~\ref{the:const}. The first is that a consistent flag-walk $W$ in a  map $\M$ whose skeleton is $q$-valent is determined by the triple $(\Phi, j, \gamma)$, where $\Phi$ is the first flag appearing in $W$, $j$ is an integer $1\leq j\leq q-1$, and $\gamma$ is one of the permutations $\alpha_j$ or $\beta_j$. For convenience we define $W(\Phi,j,\gamma)$ to be simply the consistent walk $(\Phi,\Phi^\gamma,\ldots,\Phi^{\gamma^{k-1}})$ where $k$ is the order of $\gamma$.
  
The second consequence is that since every consistent flag-walk $W$ is a $j$-hole or a $j$-Petrie path, $W$ has a well-defined reverse, namely, the reverse of the corresponding $j$-hole or $j$-Petrie path. Then, we say that $W$ is {\em symmetric} provided that there exists $g \in \Aut(\M)$ such that $g$ maps $\W$ to its reverse $\rev{W}$. We say that $W$ is {\em chiral} otherwise. 

\begin{lemma}
Let $\M$ be map of valence $q$, let $G \leq \Aut(\M)$ and let $\W = W(\Phi, j, \gamma)$ and $\W'=W(\Phi', j', \gamma')$ be two $G$-consistent flag-walks. Then, $\W$ and $\W'$ are in the same $G$-orbit if and only if $\Phi$ and $\Phi'$ are in the same $G$-orbit,  $j = j'$ and $\gamma = \gamma'$.
\end{lemma}

\begin{proof}
Let $W=(\Phi_0,\ldots,\Phi_{k-1})$ and $W'=(\Phi'_0,\ldots,\Phi'_{k'-1})$. Suppose $W$ and $W'$ are in the same $G$-orbit. Then clearly $k = k'$ and there exists a $g \in G$ such that $\Phi_i^g = \Phi'_i$ for all $i \in \{0,\ldots,k-1\}$. In particular by setting $i=0$ we see that $\Phi_0$ and $\Phi_0'$ are in the same $G$-orbit.

 Moreover, for all $i \in \{1,\ldots,k-1\}$ we have $\Phi_i = (\Phi_{i-1})^\gamma$ and $\Phi'_i = (\Phi'_{i-1})^{\gamma'}$ . Then 
$(\Phi'_i)^{\gamma'} = \Phi'_{i+1} = (\Phi_{i+1})^g = ((\Phi_{i})^\gamma)^g = ((\Phi_{i})^g)^\gamma= (\Phi'_{i})^\gamma$ for all $i \in \{1,\ldots,k-1\}$. 
In particular, $(\Phi'_0)^\gamma = (\Phi'_0)^{\gamma'}$. However $\gamma \in \{r_0r_1(r_2r_1)^{j-1}, r_0(r_2r_1)^j \}$ and $\gamma' \in \{r_0r_1(r_2r_1)^{j'-1}, r_0(r_2r_1)^{j'} \}$. We thus have four cases. 

If $\gamma = r_0r_1(r_2r_1)^{j-1}$ and $\gamma' = r_0r_1(r_2r_1)^{j'-1}$ then $(\Phi_0)^{01(21)^{j-1}} = (\Phi_0)^{01(21)^{j'-1}}$ and thus $(\Phi_0)^{01} = (\Phi_0)^{01(21)^{j'-j}}$. That is, $(r_2r_1)^{j'-j}$ fixes $(\Phi_0)^{01}$. Thus $j-j'$ is a multiple of the order of $r_2r_1$, which is $q$. This forces $j = j'$ and $\gamma =  \gamma'$. The case where $\gamma = r_0(r_2r_1)^{j}$ and $\gamma' = r_0(r_2r_1)^{j'}$ is analogous.

If $\gamma = r_0r_1(r_2r_1)^{j-1}$ and $\gamma' = r_0(r_2r_1)^{j'}$ then  $(\Phi_0)^{01(21)^{j-1-j'}} = \Phi_0^0$. Then $r_1(r_2r_1)^{j-1-j'}$ fixes $\Phi_0$, which is impossible. The same argument shows that the case where $\gamma = r_0(r_2r_1)^{j}$ and $\gamma' = r_0r_1(r_2r_1)^{j'-1}$ cannot happen. 
\end{proof}

We now have all the results and terminology necessary to characterize $G$-consistent walks in dart-transitive maps. We do so in two steps. First, in Theorem~\ref{the:orbits}, we deal with reflexible, chiral and half-reflexible maps (classes, $1$, $2$ and $2_{01}$, respectively). The characterization of consistent walks in the two remaining classes ($2_1$ and $2_0$), will follow at once from Theorem~\ref{the:orbits} and Remark~\ref{rem:petrieholes} (see Corollary~\ref{cor:orbits}).

\begin{theorem}
\label{the:orbits} 
Let $\M$ be a map of valence $q$ and let $G \leq \Aut(\M)$ be transitive on darts.
Then the following holds:
\begin{enumerate}
\item If $\M$ is $G$-reflexible, then for every $j$, the class of flag $j$-holes is a $G$-orbit of $G$-consistent flag-walks, as is the set of flag $j$-Petrie paths. 
Those are the only orbits of $G$-consistent flag-walks in $\M$ and there are $2(q-1)$ of them. 
All are symmetric.

\item If $\M$ is $G$-chiral, then for every $j$, the class of flag $j$-holes is a union of two $G$-orbits of $G$-consistent flag-walks. 
Those are the only $G$-orbits of $G$-consistent flag-walks in $\M$ and there are $2(q-1)$ of them.  
All are chiral except the lines, which are symmetric. 

\item If $\M$ is $G$-half-reflexible, then for every odd $j$, the class of flag $j$-holes is a union  of two $G$-orbits of $G$-consistent flag-walks.  
For even values of $j$, the same is true for the set of flag $j$-Petrie paths. 
Those are the only $G$-orbits of $G$-consistent flag-walks in $\M$ and there are $2(q-1)$ of them. 
The holes are symmetric, but the Petrie paths are chiral.
\end{enumerate}
\end{theorem}

\begin{proof}
We know from Theorem~\ref{the:const} that if a $G$-consistent flag-walk exists in a dart-transitive map, then it is a $j$-hole or a $j$-Petrie path for some $j$, and it is determined by a triple $(\Phi,j,\gamma)$.  
We need to show, in each of the three cases, that the indicated holes and Petrie paths are consistent and that none of the others are.  
We do this by considering the possibilities for the triple $(\Phi, j, \gamma)$. 
\begin{enumerate}
\item  If $\M$ is $G$-reflexible, then all of its flags are in the same $G$-orbit, so we may choose any flag to be $\Phi$.  
For any $1\leq j\leq q-1$, choose one of $\alpha_j$ or $\beta_j$ to be $\gamma$, so that $W=(\Phi,\Phi^\gamma,\Phi^{\gamma^2},\ldots,\Phi^{\gamma^{k-1}})$, where $k$ is the order of $\gamma$, is  $W(\Phi,j,\gamma$). The flag $\Phi^\gamma$ is in the same orbit as $\Phi$, as $\M$ is $G$-reflexible, and thus some symmetry $g$ sends $\Phi$ to $\Phi^\gamma$; that symmetry will act as a shunt for the consistent walk $W$. Indeed, for $0 \leq i \leq k-1$ we have $\Phi^{\gamma^i}=(\Phi^\gamma)^{\gamma^{i-1}} = (\Phi^g)^{\gamma^{i-1}} = (\Phi^{\gamma^{i-1}})^g$, so $W$ is shunted by $g$.

Now, in choosing the triple $(\Phi,j,\gamma)$, we had one option for the orbit of $\Phi$, $q-1$ choices for $j$ and two choices for $\gamma$, giving $2(q-1)$ choices in total for the $G$-orbits of the consistent walks. 
There are symmetries sending $\Phi$ to $\Phi_1^0$ and $\Phi_1^{0 2}$, so regardless of whether $\gamma$ is  $\alpha_j$ or $\beta_j$, there is a symmetry reversing the walk and so it is symmetric.  

\item If $\M$ is $G$-chiral, then there are two $G$-orbits of flags. 
We think of flags in one orbit as green, and in the other, red. 
For $i = 0, 1, 2$, flags $\Phi$ and $\Phi^i$ are different colours. 
We then have two options for the colour of the flag $\Phi$, and any of the $q-1$ possibilities for $j$ is feasible. 
But $\Phi^{\alpha_j}$ is the same colour as $\Phi$, while $\Phi^{\beta_j}$ is not. 
Thus we must choose $\gamma$ to be $\alpha$. 
For each of the $2(q-1)$ feasible choices for $(\Phi, j, \gamma)$, $\Phi$ and $\Phi^\gamma = \Phi^\alpha$ are the same colour, so some symmetry sends $\Phi$ to $\Phi^\gamma$, thus acting as a shunt for a consistent walk starting at $\Phi$. 
Because $\Phi^0$ is the opposite colour to $\Phi$, no symmetry will reverse the walk and so it is chiral.

\item If $\M$ is $G$-half-reflexible, then there are two $G$-orbits of flags. 
We can think of flags in one orbit as green, and in the other, red.  
As $\Phi, \Phi^0, \Phi^1$ are the same colour, and different in colour to $\Phi^2$, we can also think of the faces of $\M$ as bicoloured in red and green with adjacent faces opposite colours.  This forces the valence $q$ of $\M$ to be even.
	
Further, a product of the $r_i$'s connects flags of the same colour if and only if it has an even number of $r_2$'s. 
In particular, $\Phi$ and $\Phi^{\alpha_j}$ are the same colour (in the same orbit) exactly when $j$ is odd, while $\Phi$ and $\Phi^{\beta_j}$ are the same colour (in the same orbit) exactly when $j$ is even. 
Thus only the walks claimed in the statement of the theorem are eligible to be consistent walks. 
Conversely, for either of the two choices for the colour of $\Phi$, choosing $\gamma = \alpha$ when $j$ is odd and $\gamma = \beta$ when $j$ is even makes $\Phi$ and $\Phi^\gamma$ the same colour, forcing the shunt to exist and so all of the $2(q-1)$ choices for the triple to give orbits of consistent walks.

The claims about the walks being symmetric or chiral follow from the fact that the reverse of a $j$-hole is a $j$-hole of the same colour, while the reverse of a $j$-Petrie path of one colour is one of the other colour.
\end{enumerate}
\end{proof}
  
\begin{corollary}
\label{cor:orbits}
Let $\M$ be a dart-transitive map of valence $q$.
Then the following holds:
\begin{enumerate}
\item If $\Pe(\M)$ is $G$-chiral, then for every $j$, the class of flag $j$-Petrie paths is a union of two $G$-orbits of consistent walks. 
Those are the only $G$-orbits of consistent flag-walks in $\M$ and there are $2(q-1)$ of them.  
All are chiral except the case in which $q$ is even and $j = \frac{q}{2}$, in which case walks in that $G$-orbit are symmetric.

\item If $\Pe(\M)$ is $G$-half-reflexible, then for every odd $j$, the class of flag $j$-Petrie paths is a union  of two $G$-orbits of consistent walks.  
For even values of $j$, the same is true for the set of flag $j$-holes paths. 
Those are the only $G$-orbits of consistent flag-walks in $\M$ and there are $2(q-1)$ of them. 
The Petrie paths are symmetric, but the holes are chiral.
\end{enumerate}
\end{corollary}

\section{Types of $j$-holes: cycles, twinings, beads and bracelets}
\label{sec:cbb}

For the remainder of this paper, we will be concerned with the edge-sets of flag consistent walks in a map $\M$, and with the orbits of these sets under $\Aut(\M)$ and its dart-transitive subgroups.  Hence when we say ``$j$-hole", we mean the edge-set of a flag $j$-hole, and will use similar modifications of other terms.

Since every $G$-consistent  flag-walk is a hole of the map or its Petrie, it will be convenient to distinguish between different kinds of $j$-holes (and thus, different kinds of $G$-consistent  walks). 

Before beginning this classification, we wish to introduce notation for some configurations of edges which occur here.  The first is a {\em cycle}, the edge-set of a closed walk  having
 no  repetitions of vertices or edges. 
 Notice that a cycle may have as few as two edges.\color{Black}

Second is a {\em twining}. 
This is a closed walk  of length $2k$ with edges $[e_0, e_1,e_2, \dots, e_{2k-1}]$ such that:
\begin{enumerate}
\item it passes through $k$ distinct vertices $v_0, v_1, v_2,\dots, v_{k-1}$,
\item for all $0\leq i\leq k-1$, $e_i$ and $e_{i+k}$ have the same endvertices, namely $v_i$ and $v_{i+1}$,
\item there are numbers $d, d'$ such that at $v_i$, edges $e_i$ and $e_{i+k}$ subtend $d$ faces,  while at $v_{i+1}$, edges $e_i$ and $e_{i+k}$ subtend $d'$ faces, and
\item there is a symmetry $\mu$ of $\M$ which fixes each $v_i$ while interchanging each $e_i$ with $e_{i+k}$.
\end{enumerate}

Next is a  {\em bead}. 
We may describe a bead as an ``evenly spaced dipole". More precisely, a {\em $d$-bead} is a set of edges having the same two end vertices, say $u$ and $v$, such that at $u$ (and at $v$)  each pair of consecutive edges subtend $d$ faces. 

Finally, a {\em bracelet} is a union of $d$-beads, for some fixed $d$, which has the same vertices as some  cycle.

\begin{figure}[hhh]
\centering
\begin{subfigure}[b]{1in}
\centering
\includegraphics[height=10mm]{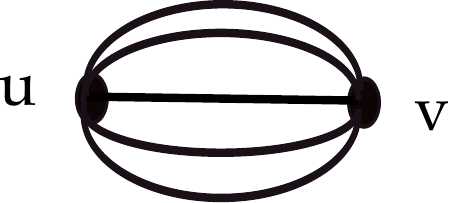}
\caption{A bead.}
\end{subfigure}
~
\begin{subfigure}[b]{2.5in}
\begin{center}
\includegraphics[height=15mm]{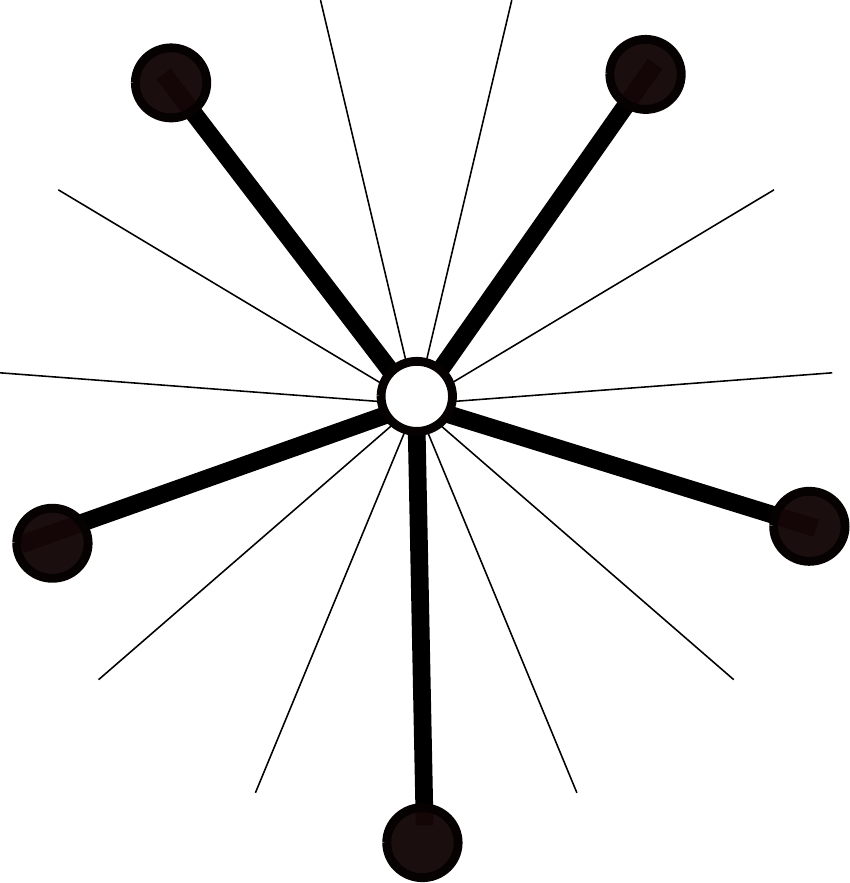}
\caption{A 3-bead at a vertex of valence 15.}
\end{center}
\end{subfigure}
~
\caption{Two views of a bead.}
\label{fig:Md}
\end{figure}

\begin{figure}[h!]
\begin{center}
\epsfig{file=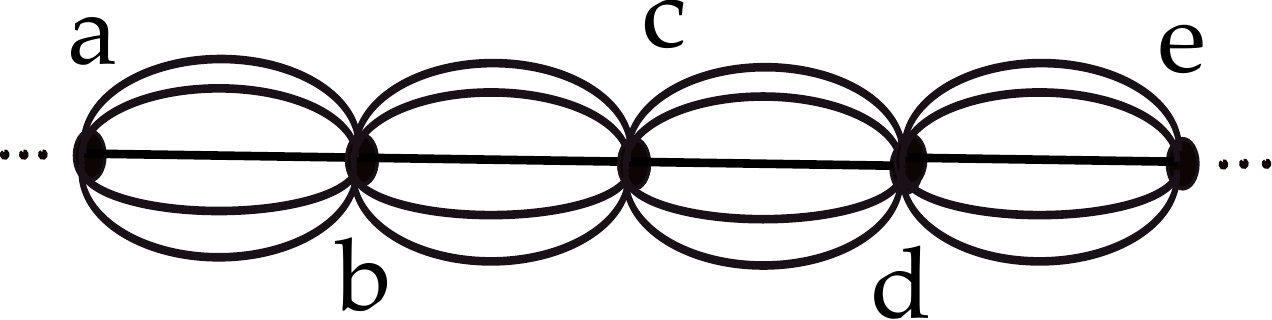,height=10mm}
\caption{Part of a bracelet.}
\label{bracelet}
\end{center}
\end{figure}

The circuit underlying a bracelet might have as few as two vertices. In this case, a bracelet on two vertices is the union of two beads, each of them on the same two vertices.  The union may or may not be a bead.

We will  consider the different ways in which a flag $j$-hole can visit an edge. We will deal with maps having only one vertex separately, in Section~\ref{sec:onevertex}. For now, we will only consider the maps that have at least two vertices, so that every edge has two distinct endvertices. 

Before we show that every $j$-hole is a twining, a cycle, a bead, or a bracelet, we offer three example maps to illustrate some possibilities.   

\subsection{Example: The map $\M = M'_{12,7}$}\label{sec:M127}
\begin{figure}[hhh]
\begin{center}
\epsfig{file=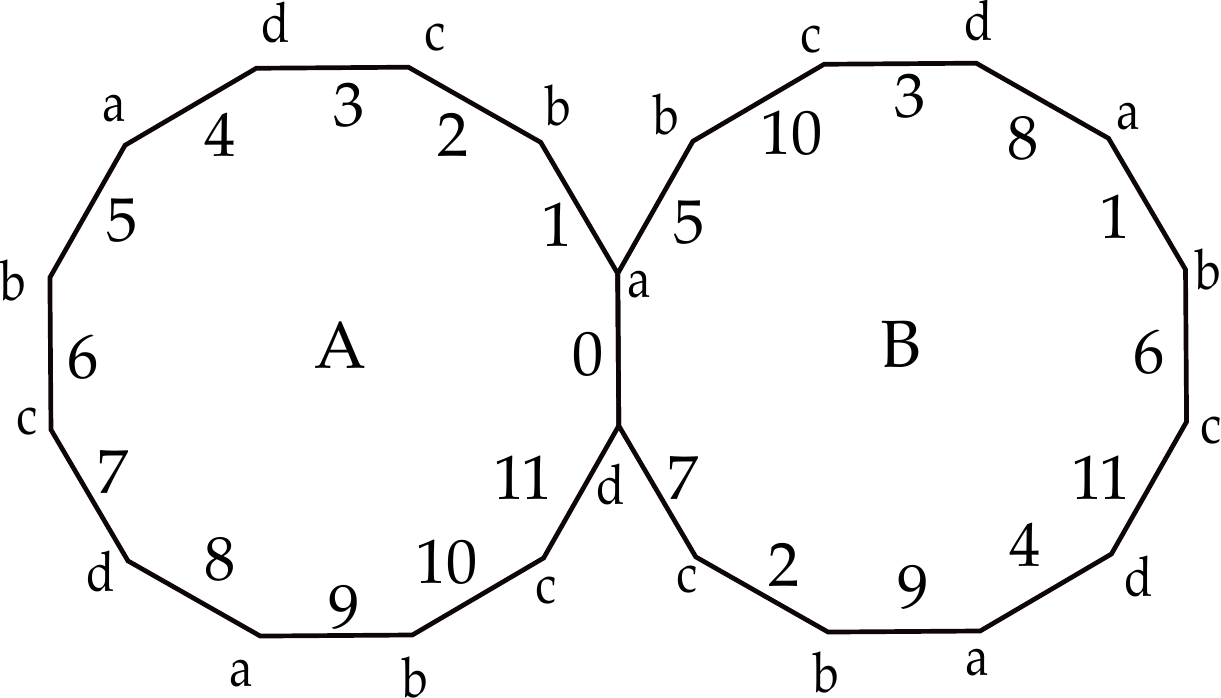,height=40mm}
\caption{A  $2$-face map.}
\label{fig:Smallmap}
\end{center}
\end{figure}

We wish to give examples showing that a $j$-hole might be a bead or a 
bracelet. One map will do for both; this is the map $\M = M'_{12,7}$, 
shown in Figure~\ref{fig:Smallmap}. For the definition of $\M_{k,i}$, see \cite{BCT}.

The map has two faces, labelled $A, B$.  It has four vertices, $a, b, c, d$, each of valence 6, and 12 edges.   It is orientable, and so lies on the orientable surface of genus 4. Its skeleton is a  4-cycle in which each edge is replaced by three edges.

\begin{figure}[hhh]
\begin{center}
\includegraphics[width=0.6\textwidth]{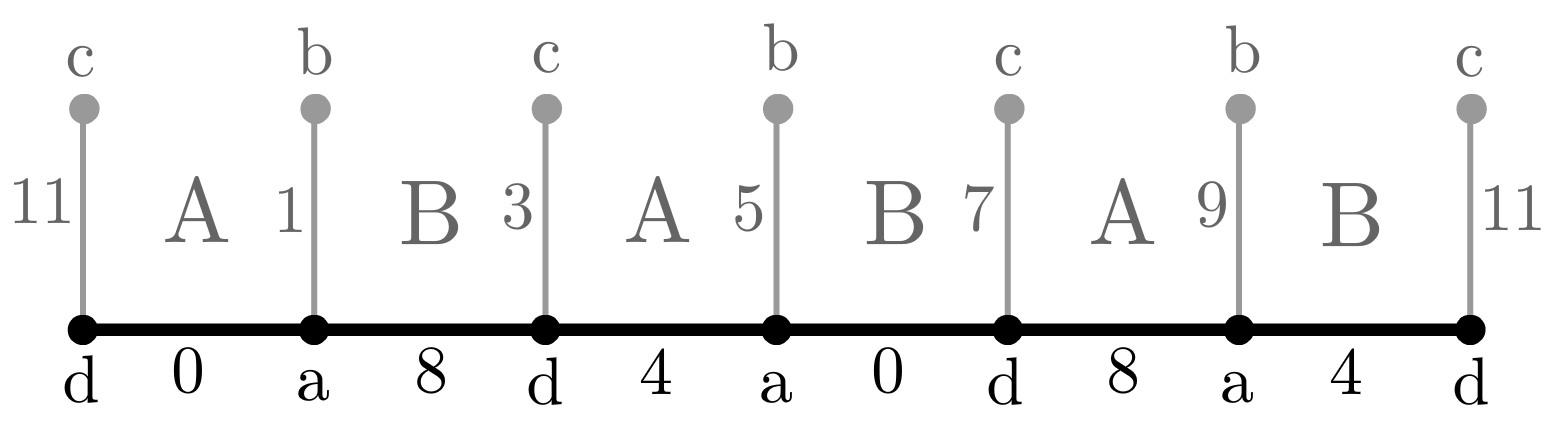}
\caption{A  second-order hole in $\M$.}
\label{fig:Small2hole}
\end{center}
\end{figure}

We show part of a 2-hole in $\M$ in Figure~\ref{fig:Small2hole}.

We see that the 2-hole alternates visiting vertices $a$ and $d$, subtending two faces  at each, and so is a bead.   The second power of the shunt along this hole fixes face $A$ and acts as a 4-step rotation there.   The third power of the shunt fixes edge 0 and acts as a $180^2$ turn there.  Thus, the hole has length 6, and its bead uses all three of the $d-a$  edges, twice each.

The reader may convince himself that a shunt along a 5-th order hole meets vertices and edges in this pattern:
$$d - 0 - a - 5 - b - 10 - c - 3 - d - 8\dots$$

It fixes face $B$ and acts as a rotation about the center of $B$, and so has order $12$.
Here, the hole visits vertices $a, b, c, d, a, b, \dots$ in that order. Thus the edges of the hole form a bracelet.

Also note that a 1-hole (i.e., a face) is a bracelet in this map.

\subsection{Example: The map $\N = \Du(M'_{12,7})$}\label{sec:DM127}

This example is provided to show that a bracelet based on a cycle of length 2 might  not be a bead.

\begin{figure}[hhh]
\begin{center}
\includegraphics[width=0.6\textwidth]{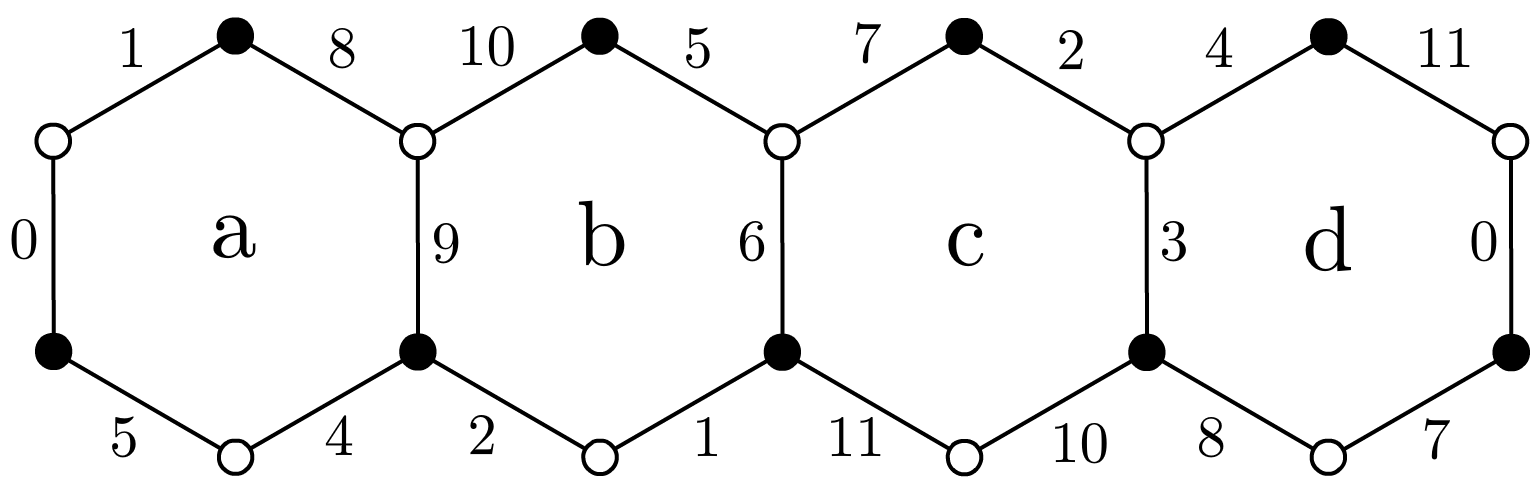}
\caption{The map $\Du(M'_{12,7})$}
\label{fig:DM127}
\end{center}
\end{figure}

The map is the dual of the previous map; $\N = \Du(\M)$.  Under the dual operator, the faces $A$ and $B$ become vertices dark and light, respectively, while the vertices $a, b, c, d$ of $\M$ become faces of $\N$.   The map $\N$ is shown in  Figure \ref{fig:DM127}.

\begin{figure}[hhh]
\begin{center}
\includegraphics[width=0.4\textwidth]{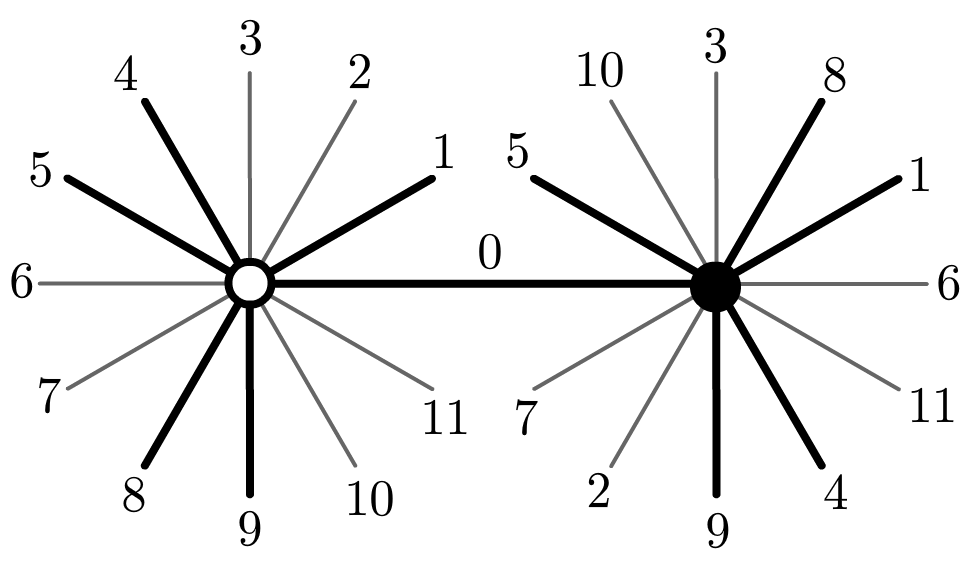}
\caption{The edges of face $a$ of $\Du(M'_{12,7})$}
\label{fig:DM127b}
\end{center}
\end{figure}

Consider face `a' in Figure \ref{fig:DM127}. It is a face and any face is a 1-hole, and thus a consistent walk.  
Its edges are 0, 1, 8, 9, 4, 5.  These are shown darkened in Figure\ref{fig:DM127b}.    Here, the edges 0,4,8 form a bead, as do 1,5,9.

 But the union of these two beads is not a bead, as the 6 edges are not equally spaced at $A$ (black dot) or $B$ (white dot).

\subsection{Cunningham's map}\label{sec:GCM} 

\begin{figure}[hhh]
\begin{center}
\epsfig{file=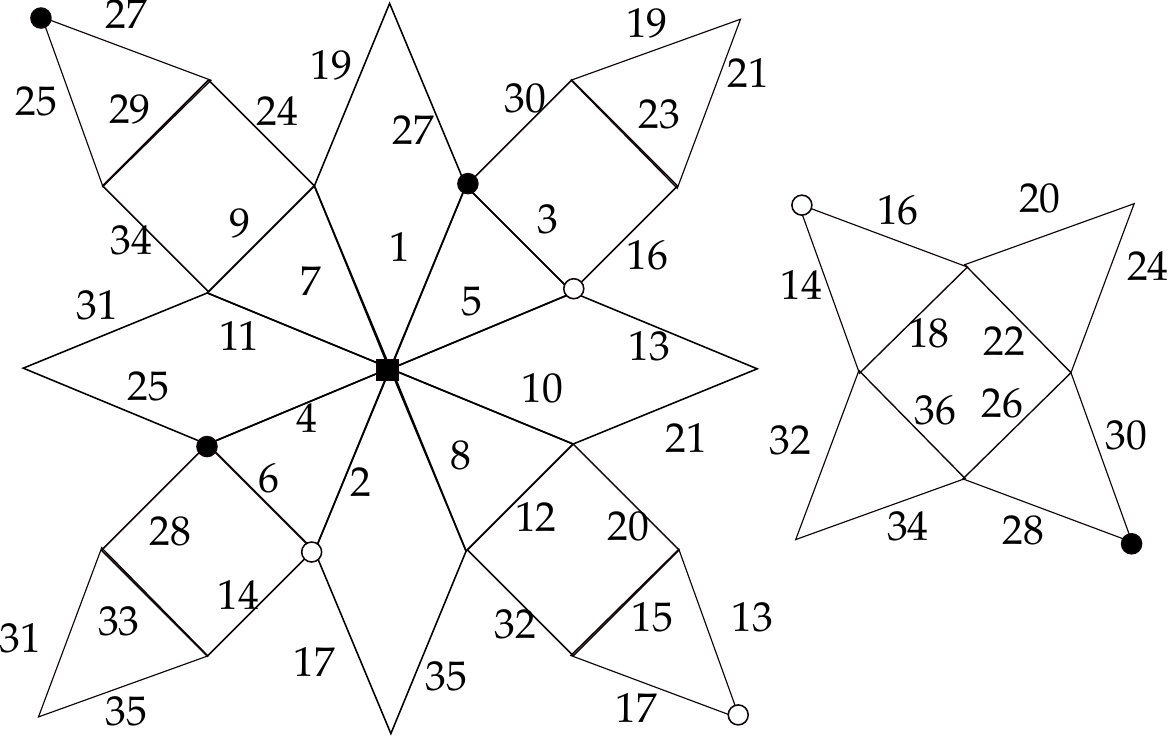,height=60mm}
\caption{A map due to Gabe Cunningham.}
\label{fig:GCMap}
\end{center}
\end{figure}

We exhibit a construction due to Gabe Cunningham \cite{GHM} 
in which a line of a map may  be a twining. The map is shown in Figure~\ref{fig:GCMap}.

\begin{figure}[hhh]
\begin{center}
\epsfig{file=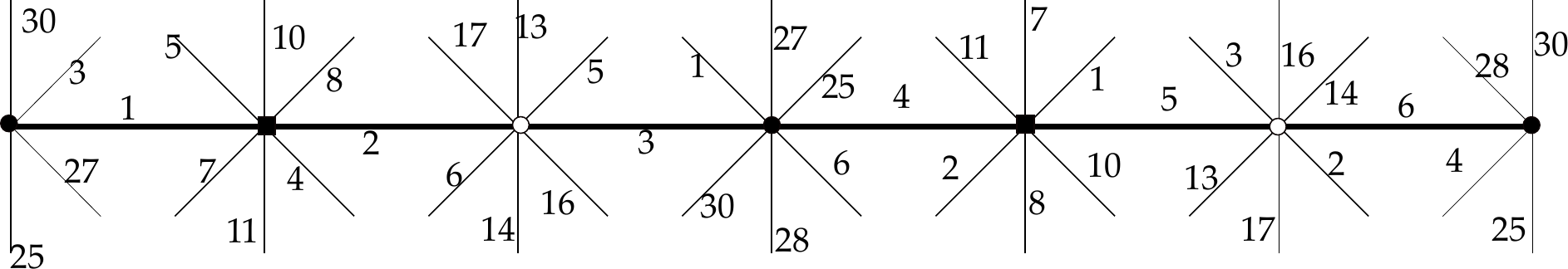,height=20mm}
\caption{A line in Cunningham's map.}
\label{fig:GCMLine}
\end{center}
\end{figure}

This map is orientable and half-reflexible.  The edges of the map are numbered so that each of $[1,2,3,4,5,6]$, $[7,8,9,10,11,12]$, $[13,14,15,16,17,18]$, $[10,20,21,22,23,24]$, $[25,26,27,28,29,30]$, $[31,32,33,34,35,36]$ is a line. The first of these is shown in Figure~\ref{fig:GCMLine}.

We can see from the figure that the line $[1,2,3,4,5,6]$ crosses itself at each of its three vertices.
The skeleton of this map is formed from $C_3\times C_3$ by doubling the edges.  The edges of each line, then, are the doubles of a 3-cycle.

\begin{theorem}\label{th:classholes}
Let $\M$ be a dart-transitive map of valence $q$ with at least two vertices. 
Then for every consistent closed flag-walk in $\M$ its edge-set
 is either a cycle, a twining, a bead, or a bracelet.
\end{theorem}
\begin{proof}
Recall from Theorem~\ref{the:const} that every consistent closed flag-walk is a $j$-hole or a $j$-Petrie path for some $j$, $0<j<q$.
As a $j$-Petrie path in $\M$ is a $j$-hole in $\Pe(\M)$, we need only concern ourselves with holes.

Let $\M$ be a dart-transitive map and let $W = (\Phi_0, \Phi_1, \dots, \Phi_{m-1})$ be a  flag $j$-hole and suppose it is shunted by some automorphism $\sigma$.   Thus, if $\Phi = \Phi_0$, then $\Phi_i = \Phi^{\sigma^i}$ for all $i$.  Now, suppose $W$ visits $\ell$ distinct edges. Then, we have four possibilities for the location of $\Phi_\ell$.

\begin{enumerate}
\item  $\Phi_\ell$ might be $\Phi$ itself. This is the simplest and most common case, and we call such a hole a {\em circuit}.   Note that, while $W$ does not visit the same edge twice, if may visit a vertex once, twice or several times. 

 If each vertex is visited once, then the circuit is a cycle.   \color{Black}   If not, let $k$ be the least positive integer such that $\Phi_0$ and $\Phi_k$ belong to the same vertex.  Let $v$ be that vertex.  We may consider two cases concerning their  relative location at $v$.

\begin{enumerate}
\item One of the flags $\Phi_0$, $\Phi_k$  might be the image of the other under the connection $(r_1r_2)^d$ for some $d$.  In that case, $\sigma^k$ acts as a $d$-step rotation about $v$.  If the $W$ visits just two vertices, its edge-set is a bead   or a pair of beads. Because a circuit on two vertices has even length, if this $j$-hole is a bead, it must have an even number of edges.  If the circuit visits more than two vertices, its edge-set is a bracelet.  The bracelet has $k$ beads, each with $\frac{\ell}{k}$ edges.\color{Black}
\item One flag might be the image of the other under the connection $(r_1r_2)^{d-1}r_1$ for some $d$.  In that case, $\sigma^k$ acts as a reflexion about some axis through $v$.  Then  $\sigma^k$ has order 2, and so $W$ has length $2k$.  Letting $e_0, e_{-1}, e_k, e_{k-1}$ be the edges containing 
$\Phi_0, \Phi_{-1}, \Phi_k, \Phi_{k-1}$, we see that at $v$, $e_0$ and $e_k$ subtend $d$ faces, while $e_0, e_{-1}$ subtend $j$ faces, as do $e_{k-1}$ and $e_k$.   This leaves $d' =q-2j\pm d$ faces to be subtended by $e_{-1}$ and $e_{k-1}$.  Because $\sigma$ is a symmetry, the corresponding edges at each vertex must subtend, $d, j$ or $q-2j\pm d$ faces.  As $\sigma^k$  fixes the vertices of the walk while interchanging the matching pairs, the $j$-hole is a twining.
 \end{enumerate}

\item  $\Phi^{\sigma^\ell}$ might be $\Phi^{0 2}$.  This case, also, is moderately common.   Let $u$ and $v$ be the (distinct) endvertices of the edge containing $\Phi$, and suppose that $ \Phi_0$ is contained in $u$, and so $\Phi_1$  is contained in $v$.   Thus, $u^\sigma = v$.    But then $\Phi_\ell = \Phi^{0 2}$ is contained in $v$ and $\Phi_{\ell+1} = \Phi_\ell^{(0 1)(2 1)^{j-1}}$ is contained in $u$, forcing $v^\sigma$ to be $u$. Thus the vertices $u$ and $v$ alternate along $W$. This forces $\ell$ to be odd.  The edges of $W$, then, consist of every $j$-th edge emanating from $u$, and similarly, every $j$-th edge emanating from $v$. That is, the edges of $W$ are edges of a dipole of $\sk(\M)$ and thus $\W$ is a $d$-bead, where $d = \GCD(q,j)$.  If the map has any $d$-beads, then every edge must belong to exactly one.

\item  $\Phi^{\sigma^\ell}$ might be $\Phi^{r_2}$, and then $\Phi^{\sigma^{2\ell}}$  must be $\Phi$ itself.  This can happen only if the map is non-orientable.   Here, the $j$-hole visits each of its edges twice, and each of its vertices at least twice.  Now, $\sigma$ sends a flag in this flag $j$-hole to the same flag that $\alpha = (r_0r_1)(r_2r_1)^{j-1}$ does.  Since $\Phi$ and $\Phi_\ell = \Phi^{r_2}$  are in $W$, $\sigma$ sends them to $\Phi^\alpha$ and $(\Phi^2)^\alpha$, respectively. But since $\sigma$ is a symmetry, we must have  $(\Phi^\alpha)^{r_2} = (\Phi^{r_2})^\alpha$.   Because $\Aut(\M)$ acts regularly on flags, this implies that $\alpha r_2 = r_2\alpha$.  Simplifying that shows that $(r_1r_2)^{2j}$ is the identity, and so $2j = q$.  Thus  such a $j$-hole  is  a line.

\item  $\Phi^{\sigma^\ell}$ might be $\Phi^{0}$. We will show that this can, in fact, not happen. To see that, first note that if $u$ and $v$ are the endvertices of the edge $e$ containing $\Phi$, and if $\Phi$ is contained in $u$, then $\Phi^0$ is contained in $v$, and so $\sigma$ first sends $u$ to $v$ and then $v$ to $u$, showing that the only vertices in $W$ are $u$ and $v$.   Similarly, the edge before $e$ in $W$ must be the same as the one after $e$.  Since $\Phi$ and $\Phi_\ell = \Phi^0$  are in $W$, $\sigma$ sends them to $\Phi^\alpha$ and $(\Phi^0)^\alpha$, respectively.  But since $\sigma$ is a symmetry, we must have  $(\Phi^\alpha)^{r_0} = (\Phi^{r_0})^\alpha$.   Because $\Aut(\M)$ acts regularly on flags, this implies that $\alpha r_0 = r_0\alpha$.  Simplifying that shows that $\alpha^2$ is the identity, and so $\Phi^{\alpha^2} = \Phi$.  This contradicts the claim that the first flag after $\Phi_0 = \Phi$ to be in $e$ is $\Phi^0$.
\end{enumerate}
\end{proof}

\begin{remark}
We see in this proof a distinction between the flag-walks tracing out a bead of odd or of even length.  We will use the phrases `odd bead' and `even bead' to indicate a bead having an odd or an even number of edges respectively. A hole which traces out an odd bead does so by visiting each edge twice, once in each direction.  An odd hole is therefore not a circuit.  A hole which traces out an even bead visits each edge once only.
\end{remark}\color{Black}


\section{Dart-transitive maps having exactly one vertex}
\label{sec:onevertex}

In this section, we will classify dart-transitive maps having exactly one vertex, determine their consistent closed walks.  We will do this by first classifying the slightly-easier-to-visualize maps having exactly one face.

Suppose that $\N$ is a   map having $n$ edges and one face. Then the face must be a $2n$-gon with edges identified in pairs. Every edge has two sides; more formally, let a {\em side} be  a pair of flags which are 0-adjacent.  If $\M = \Du(\N)$ is dart-transitive, then $\Aut(\N)$ must be  transitive on sides.

Consider a side $s = \{\Phi,\Phi^0\}$. The side of the same edge {\em opposite} to $s$ is $s^2 = \{\Phi^2,\Phi^{02}\}$. Since $\N$ has a single face, the group $\langle r_0, r_1 \rangle$ is transitive on the flags of $\N$, and so, for some integer $\ell$ we have that $\Phi^2 = \Phi^{(01)^\ell}$ or  $\Phi^2 = \Phi^{0(10)^\ell}$. In the former case, we say that $s$ and $s^2$ are non-orientably opposite; in the latter, that they are orientably opposite. We will show that the pairs of opposite sides of $\M$ are either all non-orientably opposite, or all orientably opposite.

Suppose $s$ and $s^2$ are non-orientably opposite, and consider a side $t  = \{\Psi,\Psi^0\}$. Since $\langle r_0,r_1 \rangle$ is transitive on sides, then for some positive integer $\Psi = \Phi^{(01)^k}$ and $\Psi^0=\Phi^{(01)^k0}$. Moreover, since $\Aut(\N)$ is transitive on sides, there exists an automorphism $\alpha$ mapping $s$ to $t$. There are two ways in which this could happen: either $\alpha$ maps $\Phi$ to $\Psi$ (and $\Phi^0$ to $\Psi^0$) or to $\Psi^0$. We will deal with these two cases separately.

First, suppose $\Phi^\alpha = \Psi$. We want to show that the $2$ neighbour of $\Psi$ is $\Psi^{(01)^\ell}$

$$\Psi^2= (\Phi^\alpha)^2 = (\Phi^2)^\alpha = (\Phi^{(01)^\ell})^\alpha = (\Phi^\alpha)^{(01)^\ell} = \Phi^{(01)^k(01)^\ell}=\Psi^{(01)^\ell}.$$

This shows that the opposite of $t$ is $\{\Psi^{(01)^\ell},\Psi^{(01)^\ell0}\}$, and since our choice of $t$ was arbitrary, then any two of opposite sides are non-orientably opposite.

Now, suppose $\Phi^\alpha = \Psi^0$. Then
$$ \Psi^2 = (\Phi^{\alpha 0})^2 = (\Phi^2)^{\alpha 0} = (\Phi^{(01)^\ell})^{\alpha 0} = (\Phi^\alpha)^{(01)^\ell 0} = \Psi ^{(10)^\ell},$$
but $\Psi ^{(10)^\ell} = \Psi ^{(01)^m}$ for some $m$. Thus, $t$ is non-orientably opposite to $t^2$. We conclude that the pairs of opposite sides of $\M$ are either all non-orientably opposite, or all orientably opposite.

Consider two sides which are 1-adjacent. If there is a rotation $\rho$ which takes one to the other, then that rotation must have order $2n$, and its action on the edges must be an $n$-cycle. We may number the edges so that this action is the permutation $\rho = (0\  1\ 2\ 3 \dots n-2\  n-1)$, and so that opposite sides of the face are the two sides of one edge.  Then if one pair of opposite edges are identified orientably, they all are and vice versa.  Then there are two possible maps $M_n$ and $\delta_n$, shown in Figure~\ref{fig:Md}.

\begin{figure}[hhh]
\centering
\begin{subfigure}[b]{2in}
\centering
\includegraphics[height=50mm]{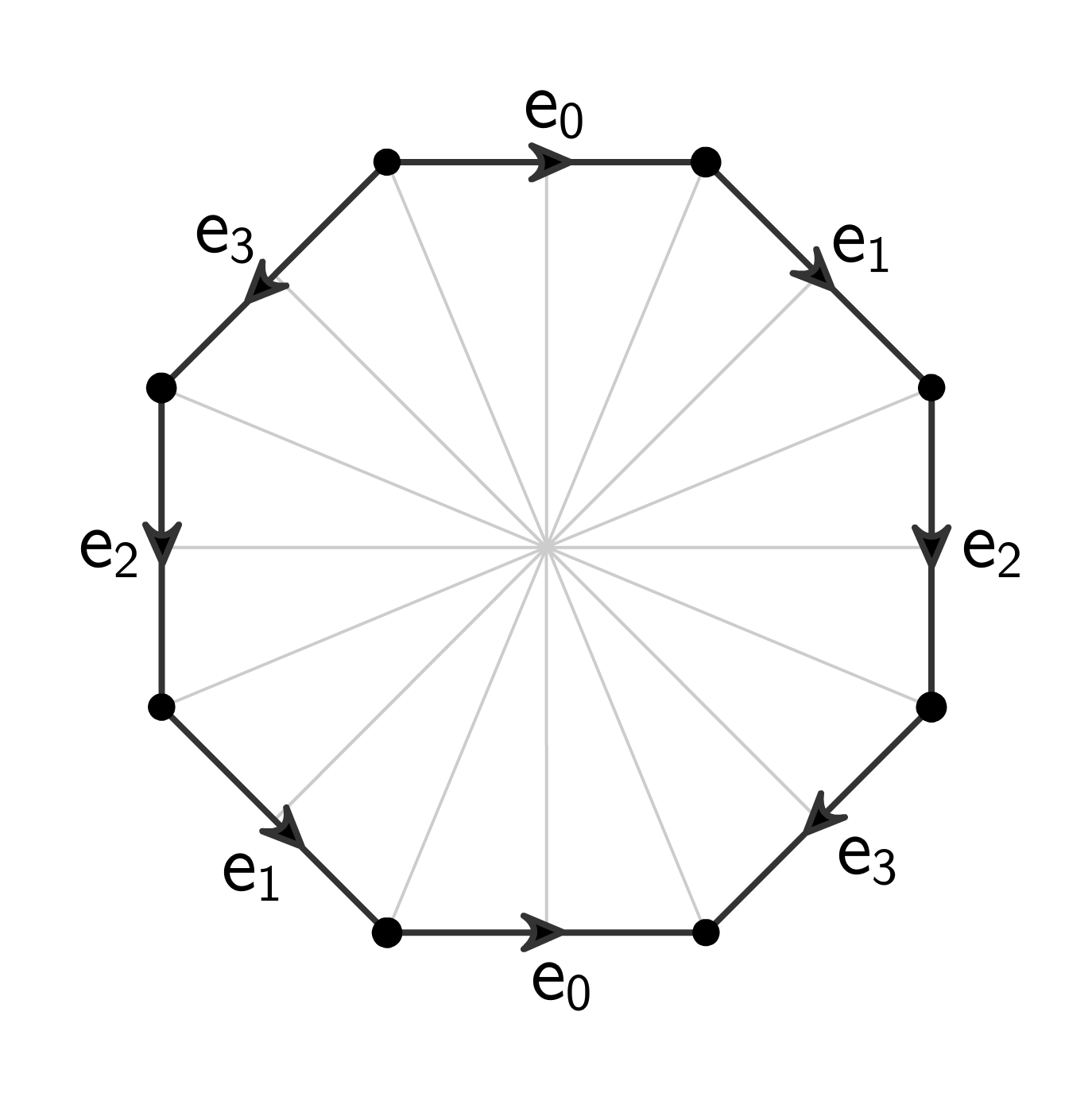}
\caption{$M_4$}
\end{subfigure}
~
\qquad \qquad
\begin{subfigure}[b]{2in}
\includegraphics[height=50mm]{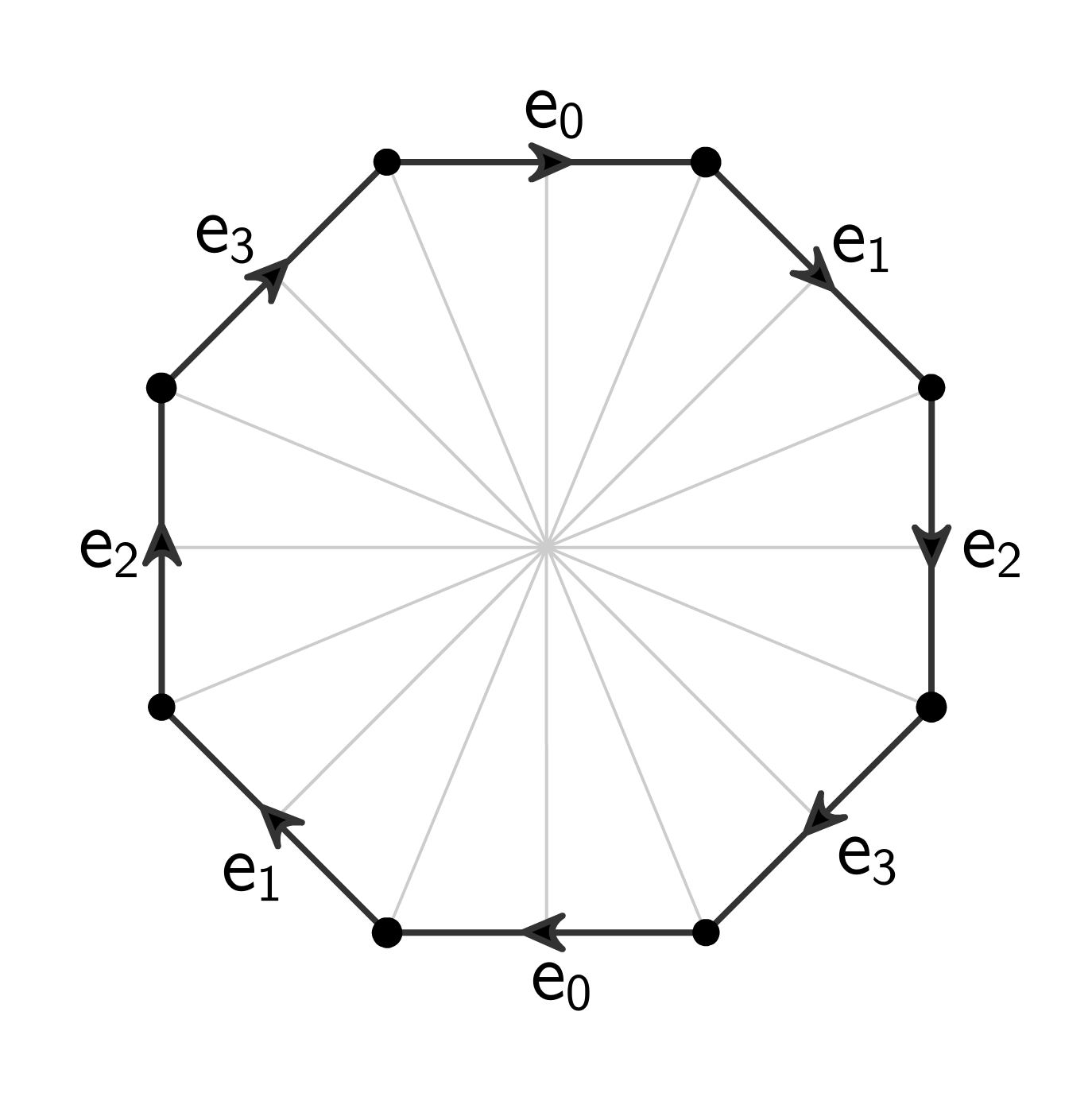}
\caption{$\delta_4$}
\end{subfigure}

\caption{Two one-face maps.}
\label{fig:Md}
\end{figure}

If there is no such rotation, then the symmetry moving one edge to the next must be a reflexion only.  Then there must be such a reflexion at each vertex, and thus there must be a symmetry $\rho'$ which acts as 2-step rotation about the face.   Then the action of $\Aut(\N)$ on the vertices has two distinct orbits, and we can colour the vertices alternately red and green.

Colour the sides alternately blue and yellow.  Then $\rho'$ permutes the sides of each colour in two $n$-cycles.  If one edge has its two sides of the same colour, then they all do. Thus, $\rho'$ acts on edges as two $\frac{n}{2}$-cycles, forcing the two sides of each edge to be opposite in the face.  Then the map is $M_n$ or $\delta_n$ again.

So, we may suppose that each edge has a blue side and a yellow side. Then  $\rho'$ cycles the edges in a single $n$-cycle.  Number the edges so that the blue sides show $0, 1, 2, \dots, n-1$ in clockwise order.
Suppose that blue 0 is between yellow $a$ and $a+1$.  With no loss of generality, assume that yellow $a$ and blue 0 share a red vertex.  The every red vertex lies between some blue $x$ and yellow $x+a$, while every green is common to blue $x$ and yellow $x+a+1$.  Then the identification of each blue $x$ with yellow $x$ is orientable.

\begin{figure}[hhh]
\begin{center}
\includegraphics[width=0.45\textwidth]{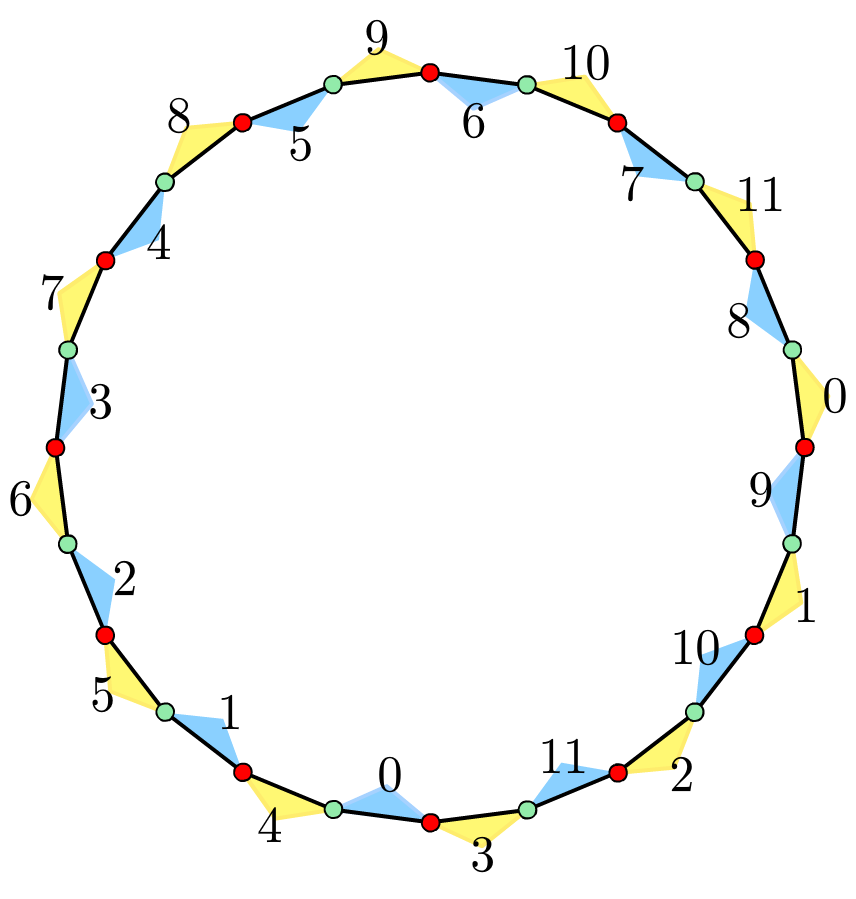}
\caption{The map $H(12,3))$.}
\label{fig:H123}
\end{center}
\end{figure}

Let $H(n,a)$ be the map so constructed. Figure~\ref{fig:H123} shows $H(12,3)$. It should be clear that $H(n, a) = H(n, -1-a).$

 If $2a+1 = n$, then $H(n,a)\cong M_n$.   
Otherwise, $H(n,a)$ has a reflexion across each axis passing through the face-center and a vertex and it has reflexions across each edge. 
 Its dual, then, is half-reflexible.   The rotation $\rho'$ is addition: $x\rightarrow x+1$, 
 and the reflexions are of the form $\mu_b: x\rightarrow b-x$ for each $b$ in $\ZZ_n$.  Then  $\Aut(\M) =\langle\rho', \mu_a\rangle$; it is isomorphic  to $D_n$.
  Around each red vertex edges appear clockwise in this order: $x, x+a, x+2a,x+3a, \dots$,
  while around each green vertex, they appear anticlockwise in order $x, x+(a+1), x+2(a+1), x+3(a+1)\dots$. 
  Thus the number of red vertices is $\GCD(n,a)$ and their valence is $|a|_n$, while green vertices, 
    $\GCD(n, a+1)$ in number, each have valence $|a+1|_n$.

\begin{figure}[hhh]
\begin{center}
\epsfig{file=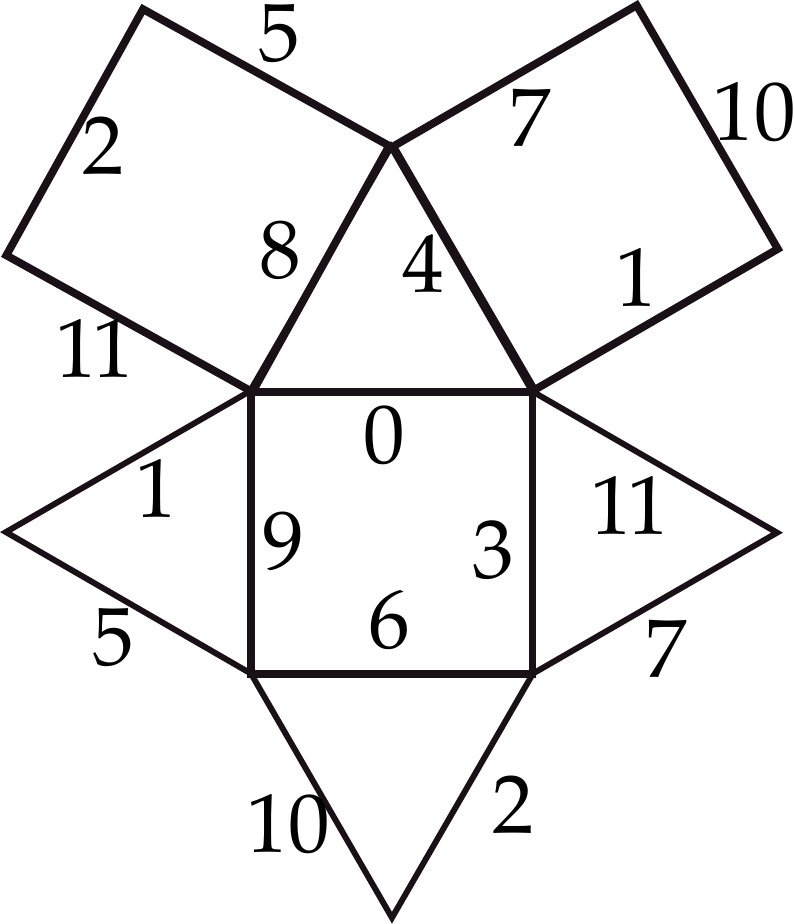,height=37mm}
\caption{The map $\Du(H(12,3))$}
\label{fig:DH123}
\end{center}
\end{figure}

For an example, consider the map $\M = \Du(H(12,3))$, shown in Figure~\ref{fig:DH123}.

\subsection{Consistent circuits and their orbits}
For $\M =   \Du(M_n)$, we determine its consistent walks.  Because $G = \M$ is reflexible, by Theorem \ref{the:orbits}(1), we know that for every $j$, every $j$-hole is consistent, as is every $j$-Petrie path.   Moreover, for each $j$, the set of all $j$-holes is an orbit as is the set of all $j$-Petrie paths.   The $j$-hole $C$ containing edge 0 visits edges $0, j, 2j, 3j, \dots$  (mod $n$).  It has length $|j+n|_{2n}$.  If we let $d$ stand for $\GCD(n,j)$, then the edges in $C$ are all multiples of $d$ (mod $n$), and the orbit of $C$ under $G$ consists of all $d$ $j$-holes.

The $j$-Petrie paths in $\M$ (other than lines) are each of length 2; each contains the edges $\{i, i+j\}$ (mod $n$) for some $i\in \ZZ_n$.  There are $n$ of these and the set of all of them is an orbit under $G$.

The map $\M$ is orientable, and so by Theorem \ref{the:orbits}(2), all of the $j$-holes are consistent   under $\Aut^+(\M)$, while none of the Petrie paths are (except lines).  The two orbits of flag-walks correspond to one orbit of edge sets.

The map $\M$ is face-bipartite if and only if $n$ is odd.  In that case, by Theorem \ref{the:orbits}(3), the $j$ holes for odd $j$, and the $j$-Petrie paths for even $j$, are consistent  under $\Aut^\#(\M)$.  The two orbits of flag-walks in each case correspond to one orbit of edge sets.
The map $\Du(\delta_n)$ is the Petrie of $\Du(M_n)$, and so has the same cycles and orbits in Petrie form.

Now we consider map $\M = \Du(H(n,a))$ for $2a+1\neq n$.  The map is half-reflexible, and so its consistent circuits are odd-order holes and even-order Petrie paths.   For the odd-order holes, if $j = 2k+1$, then one colour of $j$-holes consists of circuits of the form $i, i+a-k, i+2(a-k)\dots$, while the other colour has those of the form $i, i+a+1+k, i+2(a+1+k), \dots$.  Each colour is an orbit of  symmetric circuits (actually cycles), having lengths $|a-k|_n$ and $|a+k+1|_n$, respectively.

The  even-order Petrie paths in $\Du(H(n,a))$ all have length 2. If $j = 2k$, then each $j$-Petrie path traverses, for some $i$,  the edges $i$ and $i+k$.  These pairs are all in the same orbit 
under $\Aut(\M)$.


\color{black}


\end{document}